\documentclass[preprint,11pt]{elsarticle}

\usepackage{graphics}

\usepackage{enumitem}

\usepackage{graphicx}
\usepackage{rotating}
\usepackage{multirow}
\usepackage{epstopdf}

\usepackage{amssymb}
\usepackage{amsthm}

\usepackage{lineno}
\usepackage{fancyhdr}
\usepackage{indentfirst}
\usepackage{color}
\usepackage{newlfont}
\usepackage{textcomp}
\usepackage{amsmath}

\usepackage{latexsym}
\usepackage{mathrsfs}
\usepackage{bbold}
\usepackage{listings}
\usepackage{array}

\usepackage{algorithm,algorithmicx,algpseudocode}
\usepackage{amssymb}

\usepackage{bm}

\input epsf
\newcolumntype{T}{>{\begin{turn}{90}}l<{\end{turn}}}

\newcommand{\beq} {\begin{equation}}
\newcommand{\eeq} {\end{equation}}
\newcommand{\bit}{\begin{itemize}}
\newcommand{\eit}{\end{itemize}}
\newcommand{\bde}{\begin{description}}
\newcommand{\ede}{\end{description}}
\newcommand{\bce}{\begin{center}}
\newcommand{\ece}{\end{center}}
\newcommand{\ben} {\begin{enumerate}}
\newcommand{\een} {\end{enumerate}}
\newcommand{\bea} {\begin{eqnarray}}
\newcommand{\eea} {\end{eqnarray}}
\newcommand{\barr} {\begin{array}}
\newcommand{\earr} {\end{array}}
\newcommand{\bean} {\begin{eqnarray*}}
\newcommand{\eean} {\end{eqnarray*}}

\newcommand{\rhat}{\mathbf{\hat{r}}}
\newcommand{\mbe}{\mathbf{E}}
\newcommand{\mbx}{\mathbf{x}}

\newcommand{\mbn}{\mathbf{n}}

\newcommand{\mbf}{\mathbf{f}}

\newcommand{\mbu}{\mathbf{u}}

\newcommand{\norm}[1]{|\!| #1 |\!|}

\newcommand{\jump}[1]{\left[\!\left[ #1 \right]\!\right]_{S^*}}

\newcolumntype{C}{>{\centering\arraybackslash}p{1.5cm}}
\newcommand{\Sr}{Q}
\newcommand{\Rr}{R}

\newcommand{\Ca}{\mbox{\it Ca}}

\newcommand{\im}{{\mathrm i}}

\biboptions{sort}

\begin{document}
\bibliographystyle{abbrv}




\begin{frontmatter}
\title{A 3D boundary integral method for the electrohydrodynamics of surfactant-covered drops}

\author[kth]{Chiara Sorgentone}
\ead{sorgento@kth.se}

\author[kth]{Anna-Karin Tornberg}
\ead{akto@kth.se}

\author[nw]{Petia M. Vlahovska}
\ead{petia.vlahovska@northwestern.edu}

\address[kth]{KTH Mathematics, Linn\'e Flow Centre, 10044 Stockholm, Sweden}
	
\address[nw]{Engineering Sciences and Applied Mathematics, Northwestern University, Evanston, IL 60208, USA}

 \begin{abstract} 
We present a highly accurate numerical method based on a boundary integral formulation and the leaky dielectric model to study the dynamics of surfactant-covered drops in the presence of an applied electric field. The method can simulate interacting 3D drops (no axisymmetric simplification) in close proximity, can consider different viscosities, is adaptive in time and able to handle substantial drop deformation. For each drop global representations of the variables based on spherical harmonics expansions are used and the spectral accuracy is achieved by designing specific numerical tools: a specialized quadrature method for the singular and nearly singular integrals that appear in the formulation, a general preconditioner for the implicit treatment of the surfactant diffusion and a reparametrization procedure able to ensure a high-quality representation of the drops also under deformation. Our numerical method is validated against theoretical, numerical and experimental results available in the literature, 
as well as a new second-order theory developed for a surfactant-laden drop placed in a quadrupole electric field.
\end{abstract}

\end{frontmatter}



\section*{Introduction}
\label{Introduction}

There is currently a growing interest in the applications of electric fields to manipulate  suspensions of deformable particles. Biomedical applications span from separation and detection (for example of infected blood cells, DNA and protein molecules \cite{GASCOYNE1997240,SMITH199186,6819332}),
to selective manipulation and drug delivery (e.g. electroporation based therapies \cite{Corovic2013}). Other engineering applications are represented by mixed emulsions where a specific material needs to be isolated. A typical example is a water-in-oil emulsion: high-viscosity oils combined with asphaltenes or resins make it hard to extract the water; an electric field can then be applied for accelerating the sedimentation
process \cite{EOW2001173}. Asphaltenes, resins, waxes and similar are natural surfactants that can be found in these kind of systems, but often surfactants are also added to the emulsion to act as a demulsifier. These surface active agents are substances that modify the surface tension of the droplets and are widely used in several engineering applications to stabilize the emulsions \cite{C2SM25934F}. In this context numerical methods offer a great opportunity to better understand the physics of these systems by testing several physical situations that are often very hard to reproduce in laboratories.\\
Independently, both the influence of surfactants and electric fields on drops have been largely studied in 2D, whilst in 3D the literature is more limited. The combined effect of surfactants together with electric fields is a completely new and almost unexplored area of research in terms of numerical experiments, especially when considering multiple drops interacting.
This is due to the numerous computational challenges associated with the complex moving geometries and the multi-physics nature of the problem.\\
The effect of surfactant has been previously studied in 3D using either a boundary integral method \cite{yon_poz,BAZHLEKOV2006369,feigl}, or a diffuse-interface-method \cite{ERIKTEIGEN2011375} or again a front-tracking method \cite{MURADOGLU20082238}, but in all these cases no interactions with electric fields were considered.
The effect of electric fields on a single clean drop has been studied numerically firstly for axisymmetric configuration by Miksis \cite{Miksis}, Sherwood \cite{sherwood_1988} and later by Feng and Scott \cite{feng_scott_1996}, who included also viscosity contrast in the model. More recently, Lac and Homsy \cite{lac_homsy_2007} presented a survey of the different behaviours obtained for a wide range of resistivities and permittivities using a boundary integral method for an axisymmetric configuration and in 2015 Nganguia et al. developed a numerical scheme based on the immersed interface method for simulating the electrohydrodynamics of an axisymmetric viscous drop \cite{Nganguia2-2015}. Feng \cite{Feng2245} also found that considering the full Melcher-Taylor leaky dielectric model with charge convection can lead to enhanced prolate deformations and reduced oblate deformations as compared to the simplest model used otherwise, where unsteady terms and surface charge convection are neglected in the conservation equation for the charge density. This kind of study has been continued by Das and Saintillan who investigated theoretically \cite{das_saintillan_2017bis} (using the small deformation theory) and numerically \cite{das_saintillan_2017} (using a three-dimensional boundary element method) the effect of charge convection. \\	
Later, the interaction between a single surfactant-covered drop with an uniform electric field has been studied by Teigen and Munkejord using a level-set method in an axisymmetric, cylindrical coordinate system \cite{Teigen2010}. When considering non-uniform electric fields the available literature is even more scarce: Deshmukh and Thaokar studied theoretically and numerically the deformation, breakup and motion of a single clean drop in a quadrupole electric field \cite{Deshmukh_2012,Deshmukh_2013} and Mandal considered also the case of a surfactant covered drop, but the range of applications is strongly limited to the case of extremely weak fields due to the assumption that the drop remains spherical \cite{Mandal_2016}.\\
The interaction of multiple clean drops placed in electric fields has been studied experimentally \cite{Dommersnes2016,Varshney} but also theoretically and numerically by Baygents et al. \cite{baygents_rivette_stone_1998}, who addressed the basic two-drop problem using integral equation methods to follow the changing drop shapes and the relative motion of the drops for an axisymmetric geometry. To our knowledge there are no previous numerical studies of multiple surfactant-covered drops interacting with electric fields in a general configuration.\\
Since the numerical and experimental literature is so limited, the first step for the validation of our numerical method is to compare with theoretical results for small deformations. As previously mentioned, a theoretical approach can be applied to study the steady drop shape in the presence of a weak electric field using small perturbation analysis. The pioneer of this kind of study is Taylor who developed a first order theory for a clean drop subjected to a weak uniform electric field \cite{Taylor159}, later extended to second order by Ajayi \cite{Ajayi499}. When considering also the effect of an insoluble surfactant, first and second order small deformation theories are available for a uniform electric field \cite{Nganguia2013,Ha1995}, and a first order theory has been developed by Mandal et al. \cite{Mandal_2016} for the case of a quadrupole field. In the present work we will also present an extension to the second order to enable comparisons of our numerical results with a more accurate theory also for quadrupole fields. A recent review that summarizes experimental and theoretical studies in the area of fluid particles (drops and vesicles) in electric fields with a focus on the transient  dynamics and  extreme  deformations is also available \cite{Petia_review}. 
\\
In the presence of an electric field, the interfacial force will depend not only on the geometry and surfactant concentration, but also on the electric stresses acting on the drop surfaces which need to be computed by solving the corresponding partial differential equation; for doing this, we will extend the numerical method presented in \cite{sorgentone_tornberg} to include the effect of the electric field. \\
The paper is organized as follows: in section \ref{sc:mathematical_formulation} we introduce the model and the boundary integral formulation used; in section \ref{sc:num_method} we discuss the numerical method used to solve the system for the evolution of the drops and in section \ref{sc:spt} a summary of the small perturbation theories available in literature for clean and surfactant covered drops in electric fields is given with an extension to the second order theory for surfactant covered drops in a quadrupole field. The numerical method is validated in section \ref{sc:num_exp} againts the small perturbation theories, the spheroidal model \cite{Nganguia2013}, numerical experiments \cite{Teigen2010} and experimental results \cite{HA1998195}. In this section we also show the spectral convergence of our method and a three drops interaction numerical experiment. Conclusions and future work are discussed in the last section.

\section{Mathematical formulation}
\label{sc:mathematical_formulation}

We consider $N$ drops suspended in an ambient fluid. The Stokes equations read 
\begin{equation}
\label{eq:diff_eq1}
\begin{cases}
-\mu_i \Delta \textbf{u}(\textbf{x})+\nabla P(\textbf{x}) =0 \\ 
\nabla \cdot \textbf{u}(\textbf{x})=0 
\end{cases}
\end{equation}
for every $\textbf{x}$ inside the $i$-th drop ($i=1,\dots,N$) or in the exterior region ($i=0$), where $\textbf{u}$ is the fluid velocity, $P$ is the pressure and $\mu_i$ is the viscosity. The fluid motion is coupled to the interface motion via the kinematic boundary condition,
\beq \dot{\mbx} = \mbu(\mbx), \text{ for all } \mbx \in S^*, \label{eq:dxdtugualeu} \eeq
where $\mbx$ is the membrane position and $S^*=\bigcup_i S_i$ denote the union of all drop surfaces. The permittivity $\epsilon$ and the conductivity $\sigma$ are discontinuous across the interface. A stress balance at the interface establishes the flow and electric field interaction:
\beq  \jump{\mbn\cdot(\Sigma^{el} + \Sigma^{hd})} = \mbf 
\label{eq:stress_bal}
\eeq
where $\jump{\cdot}$ denotes the jump across the interface (e.g., $\left[\!\left[ \sigma \right]\!\right]_{S_i}= \sigma_0 - \sigma_i$), $\mbn$ is the outward pointing unit normal, $\Sigma^{el}$ is the electric stress, $\Sigma^{hd}$ is the hydrodynamic stress.\\
Denoting by $\gamma$ the interfacial tension, the interfacial force is defined by:
\beq
\label{eq:interfacial_force}
 \mbf = 2 \gamma(\mbx)H(\mbx)\mbn(\mbx) - \nabla_S \gamma,
\eeq
where $H=\frac{1}{2} \nabla_S \cdot \mbn$ denotes the mean curvature and $\nabla_S=(I-\mbn \mbn) \nabla$ is the surface gradient. For a clean drop, the surface tension coefficient $\gamma(\mbx)$ will be constant, and the second term in (\ref{eq:interfacial_force}), the so-called Marangoni force, will vanish.\\
The electric stress $\Sigma^{el}$ is given by the Maxwell stress tensor, defined as,
\beq
\Sigma^{el} =  \tilde{\epsilon}\epsilon_0( \mbe \otimes \mbe - \frac{1}{2} \norm{\mbe}^2 \,  \mathbf{I})
\label{maxwell} \eeq
where $\tilde{\epsilon}$ is the permittivity of the vacuum; the hydrodynamic stress tensor is given by 
\beq
\Sigma^{hd}=-PI+\mu(\nabla \mbu+\nabla \mbu^T).
\label{Cauchystresstensor} \eeq

To solve for the drop evolution under flow and electric fields, we use a boundary integral equation (BIE) method. 
Letting $a$ and $\mu_0$ be, respectively, the radius of an undeformed spherical droplet and the viscosity of the ambient fluid and defining $E_\infty=\sqrt{\frac{Ca_E \gamma_{eq}}{a\epsilon_0 \tilde{\epsilon}}}$, where $Ca_E$ is the electric capillary number, we can
non-dimensionalize by taking the characteristic length to be $a$, the characteristic electric field to be $E_\infty$, the
characteristic velocity to be $U=\frac{\epsilon_o \tilde{\epsilon} E_\infty^2a}{\mu_0}$ and the
characteristic time $T=\frac{a}{U}$ \cite{baygents_rivette_stone_1998}, where $\gamma_{eq}$ is the equilibrium surface tension. Using the
corresponding non-dimensionalized variables, equations
(\ref{eq:diff_eq1}) and (\ref{eq:stress_bal}) can be reformulated as a
boundary integral equation \cite{pozrikidis_1992}. For all $\textbf{x}_0 \in S_i \ (i=1,\dots,N)$,
 \begin{equation}
 \label{eq:main_eq}
 \begin{aligned}
 (\lambda_i+1)\textbf{u}(\textbf{x}_0)=&-\sum_{j=1}^N  \bigg( \frac{1}{4\pi}\int_{S_j} \bigg(\frac{\textbf{f}(\textbf{x})}{Ca_E}-\textbf{f}^E(\textbf{x})\bigg)\cdot G(\textbf{x}_0,\textbf{x})dS(\textbf{x})\bigg) \\
 &+\sum_{j=1}^N \bigg(\frac{\lambda_i-1}{4\pi}\int_{S_j} \textbf{u}(\textbf{x})\cdot T(\textbf{x}_0,\textbf{x})\cdot \textbf{n}(\textbf{x}) dS(\textbf{x})\bigg),
 \end{aligned}
\end{equation}
where $\mbf$ was defined in (\ref{eq:interfacial_force}), $\mbf^E$ is the electric force on the interface $\mbf^E = \jump{\mbn\cdot\Sigma^{el}}$ and $\lambda_i=\frac{\mu_i}{\mu_0}$ denotes the viscosity contrast of the $i$-th drop.
The tensors $G$ and $T$ are the Stokeslet and the Stresslet, 
\beq
\label{eq:GandT}
G(\textbf{x}_0,\textbf{x})=I/r+\hat{\textbf{x}}\hat{\textbf{x}}/r^3, \ \ \ T(\textbf{x}_0,\textbf{x})=-6\hat{\textbf{x}}\hat{\textbf{x}}\hat{\textbf{x}}/r^5,
\eeq
with $\hat{\textbf{x}}=\textbf{x}_0-\textbf{x}$ and
$r=|\hat{\textbf{x}}|$. Note that we are assuming there is no far-field velocity so that the motion of the drops is driven only by the electric field, but the extension considering also a background flow is straightforward \cite{pozrikidis_1992,sorgentone_tornberg}. \\
The evolution of the drops is influenced by the effect of the electric field that we will study using the {\em leaky--dielectric} model, in which the electric charges are assumed to be present only at the interface and not in the bulk. The boundary value problem for the electric field can be written as:  
\begin{subequations}
\begin{align}
-\nabla \cdot \mathbf{E} = 0   \quad\text{in}\quad \mathbb{R}^2\setminus S^* &\\
\jump{\sigma E_n} = 0 &  \label{sub-eq-1:9b} \\ 
\mathbf{E}(\mbx) \rightarrow \mbe_\infty(\mbx) \quad \text{as} \quad \norm{\mbx} \rightarrow \infty & 
\end{align} \label{eqn:potential}
\end{subequations}
where $E_n=\mathbf{E}\cdot \mathbf{n}$ and $\mathbf{E}_t=(I-\mbn \mbn) \cdot \mathbf{E}=\mathbf{E}-E_n\mathbf{n}$ are respectively the normal and tangential components of $\mathbf{E}$. In the present paper we consider that the surfactant-covered drops are subjected to the electric field $\textbf{E}_\infty$ applied far away from the drop due to an electric potential of the form 
$$\phi_\infty=-[E_u \bar{r} P_1(\eta)+E_q \bar{r}^2 P_2(\eta)],$$
where $\bar{r}$ is the radial coordinate and $E_u$ and $E_q$ represent the strength of uniform and quadrupole component respectively. $P_n(\eta)$ is the Legendre polynomial of degree $n$ with argument $\eta=cos(\theta)$, where $\theta$ is the zenith angle.
Note that the boundary condition \eqref{sub-eq-1:9b} is obtained by neglecting unsteady terms and surface charge convection to the conservation equation for the charge density $q=\jump{\epsilon E_n}$ \cite{das_saintillan_2017}.\\
For simplicity we will assume the viscosity ratio $\lambda=\frac{\mu_i}{\mu_0}$, the conductivity ratio $R=\frac{\sigma_i}{\sigma_0}$ and the permittivity ratio $Q=\frac{\epsilon_i}{\epsilon_0}$ to be the same for all the drops, $i=1,\dots,N$. Using these definitions, eq. (\ref{eqn:potential}b) implies 
\beq E_n^0=RE_n^i. \label{eq:EneEni}\eeq
We will henceforth omit the superscript $0$ for the normal component of the electric field.\\
The electric force on each interface $i$ can be written in terms of $E_n$ and $\textbf{E}_t$ \cite{lac_homsy_2007}:
\beq
\label{eq:el_force}
\mbf^E = \jump{\mbn\cdot\Sigma^{el}}
=\frac{\tilde{\epsilon}\epsilon_0}{2}[(1-\frac{Q}{R^2})E_n^2-(1-Q)E_t^2] \mathbf{n}+\tilde{\epsilon}\epsilon_0(1-\frac{Q}{R})E_n\mathbf{E}_t.
\eeq
The electric force on the membrane needs to be computed by solving (\ref{eqn:potential}) for a given drop shape. Since (\ref{eqn:potential}) is a linear partial differential equation, similar to the fluid problem, we can recast it using a boundary integral formulation \cite{lac_homsy_2007,baygents_rivette_stone_1998}:
\begin{equation} 
\label{eq:BIE01}
\begin{split}
\mathbf{E}_\infty(\mbx_0) -\sum_{j=1}^N \int_{S_j} \frac{\hat{\mbx}}{4\pi r^3} {\left[\!\left[ E_n(\mbx) \right]\!\right]}dS(\mbx)= \begin{cases} 
\mathbf{E}^i &\mbox{if } \mbx_0 $ inside $ S^*, \\ 
\frac{1}{2}[\mathbf{E}^0+\mathbf{E}^i] &\mbox{if } \mbx_0 \in S^*, \\ 
\mathbf{E}^e &\mbox{if } \mbx_0 $ outside $ S^*. \\ 
\end{cases}
\end{split}
\end{equation}
Eq. (\ref{eq:BIE01}) exactly satisfies the far-field condition (\ref{eqn:potential}c) and gives an integral equation for $E_n$ by taking its inner product with the normal vector and using (\ref{eq:EneEni}):
\beq
\label{eq:E_n}
\frac{R}{R+1}\mathbf{E}_\infty \cdot \mbn(\mbx_0)+\frac{R-1}{R+1} \sum_{j=1}^N \int_{S_j} \frac{\hat{\mbx} \cdot \mbn}{4\pi r^3} E_n(\mbx)dS(\mbx)=\frac{1}{2}E_n(\mbx_0).
\eeq
The tangential component of the electric field is given by
\beq
\label{eq:E_t}
\mathbf{E}_t=\frac{\mathbf{E}^0+\mathbf{E}^i}{2}-\frac{1+R}{2R}E_n \mbn.
\eeq
The presence of an insoluble surfactant on the drop surface will affect the interfacial tension $\gamma$. It is related with the surfactant concentration by the equation of state. Different equations of state can be used \cite{Pawar}, as the Langmuir equation of state which is given, in dimensionless form, by
\beq
\label{eq:langmuir}
\gamma(\Gamma)=\gamma_0(1+\beta \ln(1-x_s\Gamma))
\eeq
where $\gamma_0$ is the surface tension of a clean drop, $\beta$ is the elasticity number and $x_s$ is the surface coverage, $0 \leq x_s \leq 1$, or the linear equation of state:
\beq
\label{eq:lin_eq_state}
\gamma(\Gamma)=1+\tilde{\beta} (1-\Gamma)
\eeq
where $\tilde{\beta}=\frac{\gamma_0-\gamma_{eq}}{\gamma_{eq}}=Ca_E Ma$, $Ma$ denotes the Marangoni number $Ma=\frac{\tilde{\beta}}{Ca_E}$, $\Ca_E=Ca_E^0(1+\tilde{\beta})$ wit $Ca_E^0$ being the capillary number based on the clean drop \cite{Nganguia2013,vlahovska_blawzdziewicz_loewenberg_2009}.\\
The equation governing the evolution of the surfactant concentration is a convective-diffusion equation which can be derived stating the conservation of surfactant mass; it is given in dimensionless terms by \cite{sorgentone_tornberg,Stone90}:
\beq
\label{eq:surfactant0}
\frac{\partial \Gamma}{\partial t}+\nabla_S \cdot (\Gamma \textbf{u}_t) -\frac{1}{Pe}\nabla_S^2\Gamma+2H(\textbf{x})\Gamma (\textbf{u}\cdot \textbf{n})=0,
\eeq
where $Pe$ is the P\'eclet number.  
\\
The steps involved in evolving the drop can now be summarized. Given $\textbf{E}_\infty$, the drop shape and the surfactant concentration:
  \begin{enumerate}[label=\roman{*}., ref=(\roman{*})]
\item Solve the system (\ref{eq:E_n}) for the unknow $E_n$;
\item Evaluate the integral (\ref{eq:BIE01}) to compute $\frac{\mathbf{E}^0+\mathbf{E}^i}{2}$;
\item Compute the tangential component of the electric field $\textbf{E}_t$ using eq. (\ref{eq:E_t});
\item Compute the electric force on the interface with eq. (\ref{eq:el_force});
\item Compute the surface tension using eq. (\ref{eq:langmuir}) and then the interfacial force using eq. (\ref{eq:interfacial_force});
\item Substituting both the interfacial and the electric force into \eqref{eq:main_eq}, solve for updated velocity;
\item Compute the new interface position, eq. (\ref{eq:dxdtugualeu});
\item Update the surfactant concentration, eq. (\ref{eq:surfactant0}).
\end{enumerate}

\section{Numerical Method} \label{sc:num_method}
In this section we will describe the numerical method that we will use for simulating 3D drops covered by an insoluble surfactant placed in an electric field. The method is an extension of our previous work \cite{sorgentone_tornberg}, where surfactant-covered drops but not electric fields were considered. The presence of an electric field requires the solution of an extra integral equation (step i. in section \ref{sc:mathematical_formulation}) and the evaluation of the mean electric field (step ii.). We will follow the general framework introduced in \cite{sorgentone_tornberg}, where a more detailed explanation of each specific numerical tool can be found.\\
Assuming the surface $S_i$ of each drop to be smooth and of
spherical topology, spherical harmonics expansions will be used to
represent the surfactant concentration $\Gamma$ and each component of the position vector $\mbx$ and the electric field $\textbf{E}$. The order of truncated spherical harmonic expansion will be denoted by $p$, and all the previous variable will be represented by
\begin{equation}
\label{eq:spharm_ex}
y_p(\theta,\phi)= \sum_{n=0}^p \sum_{m=-n}^n y_n^m Y_{nm}(\theta,\phi),
\end{equation}
where $\phi$ and $\theta$ represents a Gaussian spherical grid: in the longitudinal direction we discretize $\phi \in [0,2\pi)$ using $2p+2$ equidistant points $\left\{ \phi_i=\frac{\pi i}{p+1} \right \}_{i=0}^{2p+1}$, while in the non-periodic direction $\theta \in [0,\pi]$ we set $\left\{ \theta_j=\cos^{-1}(t_j)\right \}_{j=0}^{p}$, where $t_j \in [-1,1]$ are the Gauss-Legendre quadrature nodes with corresponding Gaussian quadrature weights $w_j^G$. The normalized scalar spherical harmonic function of degree $n$ and order $m$ is given by
\begin{equation}
\label{eq:normalized_spharm}
Y_{nm}(\theta,\phi)=\sqrt{\frac{2n+1}{(4\pi)}\frac{(n-m)!}{(n+m)!}}\tilde{Y}_{nm}, \ \ \ \text{with } \tilde{Y}_{nm}=P_n^{m}(\cos(\theta))e^{im\phi}, \ \ \ -n\leq m \leq n,
\end{equation}
where the associated Legendre function is defined by
\begin{equation}
P_n^{m}(x)=(-1)^{m}(1-x^2)^{m/2}\frac{d^m}{dx^m}P_n(x), \ \ \ m \geq 0,
\end{equation}
and
\begin{equation}
\label{eq:pneg}
P_n^{-m}(x)=(-1)^{m}\frac{(n-m)!}{(n+m)!}P_n^m(x), \ \ \ m \geq 0.
\end{equation}
$P_n$ is the Legendre polynomial of degree $n$
\begin{equation}
P_n(x)=(2^nn!)^{-1}\frac{d^n}{dx^n}[(x^2-1)^n].
\end{equation}
\\
The whole machinery for computing the evolution of the drops is solved using a Galerkin formulation as suggested in \cite{rahimian_veerapaneni}, where a similar approach is used to simulate the evolution of vesicles.\\
Let $(\cdot,\cdot)$ denote the inner product
$(y,Y_{nm})=\int_\mathbb{S^2} y \overline{Y_{nm}}ds$, then the Galerkin
method seeks the solution to the original system by
\begin{equation}
\label{eq:BIM_galerkin}
(\lambda+1)(\textbf{u},Y_{nm})= \sum_{j=1}^N [ (\textit{S}_j[\textbf{f}],Y_{nm})+((\lambda-1)\textit{D}_j[\textbf{u}],Y_{nm})],
\end{equation}
\begin{equation}
(\frac{\partial \textbf{x}}{\partial t},Y_{nm}) = (\textbf{u},Y_{nm}),
\end{equation}
\begin{equation}
\label{eq:galerkin_surfactant}
(\frac{\partial \Gamma}{\partial t},Y_{nm}) = -(\nabla_\gamma \cdot \Gamma \textbf{u}_\gamma, Y_{nm}) +\frac{1}{Pe}(\nabla_\gamma^2\Gamma,Y_{nm})-(2H\Gamma \textbf{u}\cdot \textbf{n},Y_{nm}),
\end{equation}
where the electric stress that contributes to the forcing term in (\ref{eq:BIM_galerkin}) is also computed using a Galerkin formulation. To be more precise, for the normal component of the electric field we solve the following system:
\begin{equation}
\label{eq:galerkin_En}
\frac{1}{2}(E_n,Y_{nm}) = \frac{R}{R+1}(\textbf{E}_\infty \cdot \textbf{n}, Y_{nm}) +\frac{R-1}{R+1}\sum_{j=1}^N(\textit{L}_j[E_n],Y_{nm}).
\end{equation}
The operators $\textit{S}, \textit{D}$ denote the single and double layer Stokesian operators and $\textit{L}$ is the Laplace double layer operator:
\beq
\label{eq:op_stokes}
\textit{S}_j[\mbf](\mbx_0)=\frac{1}{4\pi} \int_{S_j}  \mbf(\mbx) \cdot G(\mbx_0,\mbx)dS(\mbx), \ \ \textit{D}_j[\mbf](\mbx_0)=\frac{1}{4\pi} \int_{S_j}  \mbu(\mbx) \cdot T(\mbx_0,\mbx)\cdot  \mbn(\mbx)dS(\mbx)
\eeq
\beq
\label{eq:dl_op_laplace}
\textit{L}_j[E_n](\mbx_0)=\frac{1}{4\pi} \int_{S_j}  E_n(\mbx)  \mathcal{L}_l(\mbx_0,\mbx)\cdot  \mbn(\mbx)dS(\mbx),
\eeq
where $G$ and $T$ have been defined in (\ref{eq:GandT}) and $\mathcal{L}_l(\mbx_0,\mbx)=\frac{\hat{\mbx}}{r^3}$.
This approach results to be cheaper compared to a standard approach as observed already in \cite{rahimian_veerapaneni,sorgentone_tornberg} thanks to the orthogonality and symmetry properties of the spherical harmonics \cite{Atkinson}. Indeed the systems to solve (for computing the velocity field - eq. (\ref{eq:BIM_galerkin}) -, the electric field - eq. (\ref{eq:galerkin_En}) - and the surfactant concentration - eq. (\ref{eq:galerkin_surfactant}) - ) reduce to almost half of the original size: $(p+1)(p+2)$ for the spherical coefficients instead of $2(p+1)^2$ for the physical variables.\\
The derivatives of any physical variable are computed using analytical expressions for the spherical harmonic coefficients (see Appendix A in \cite{sorgentone_tornberg}). An upsampling and de-aliasing procedure is needed to accurately compute geometrical quantities involving non-linear manipulations and for this reason we will use the adaptive algorithm proposed by Rahimian et al. \cite{rahimian_veerapaneni} which is based on the decay of the spectrum of the mean curvature.
\subsection{Regular, singular and nearly-singular integration}
\label{sc:quadrature}
The quadrature rule for regular integrals is defined as follows:
\begin{equation}
\label{eq:reg_quadrature}
\int_{S} y d\gamma \approx \sum_{j=0}^p \sum_{k=0}^{2p+1} w_{j} y(\theta_j,\phi_k)W(\theta_j,\phi_k),
\end{equation}
where $w_{j}=\frac{\pi}{p}\frac{w_j^G}{ sin(\theta_j)}$ and $W(\theta_j,\phi_k)$ is the infinitesimal area element of the surface $S$.\\
When computing the integrals in (\ref{eq:main_eq}), (\ref{eq:BIE01}), and (\ref{eq:E_n}) we need a special treatment for the case
$\textbf{x}_0=\textbf{x}$ (where we have a singularity), and for the case $lim_{\textbf{x}_0\rightarrow \textbf{x}}$ (where we have a so-called nearly-singular integral). These two situations have a very different nature: in the first case it is an analytical problem while in the second case it is a purely numerical issue. For this reason we will treat the two situations separately. 

For the singular case, we will make use of the fact that the spherical harmonics are eigenfunctions of the Laplace operator on the sphere \cite{ganesh}. This property can be extended to the case of the single and double layer kernels using symmetries and smoothness properties to build modified weights that ensures exponential convergence \cite{rahimian_veerapaneni,veerapaneni2011}. The modified weights are given by
\beq
\label{eq:mod_weights}
w_j^{mod}=w_j \sum_{n=0}^p 2sin(\theta_j/2)P_n(cos(\theta_j)
\eeq
and they are made for the target point $\textbf{x}_0$ to be the north pole of the parametrization. Hence, for
a general target point different from the north pole, the coordinate system will be rotated such that the target point becomes the north pole and then the special quadrature rule can be applied. Analytical expressions for rotating the coefficients of the spherical harmonic expansion are available \cite{sorgentone_tornberg} and a fast algorithm for spherical grid
rotation with application to singular quadrature is given by Gimbutas and Veerapaneni \cite{gimbutas} where the computational complexity for all targets on the surface is $O(p^4\log{p})$.\\
For the nearly-singular case we will use an approach based on interpolation: the idea is that we can accurately compute the on-surface integrals using the modified weights in (\ref{eq:mod_weights}) and, using the regular quadrature (\ref{eq:reg_quadrature}) with a reasonable upsampling rate, also the integrals where the target point is far enough from the surface. We can do that in a number of aligned points and then perform a 1D Lagrange interpolation at the original target point. This procedure was firstly proposed by Ying et al. \cite{Ying_Biros} and then optimized and tailored to our specific setting in \cite{sorgentone_tornberg}. \\

\subsection{Time-stepping and reparametrization}
\label{sc:time_and_rep}
We recall that the system of ODEs (drop/surfactant) we are considering is given by
 \begin{equation}
 \label{eq:system1fefi}
\begin{cases}
\frac{d\textbf{x}}{dt}=\textbf{u}(\textbf{x},\sigma), \\
\frac{d\Gamma}{dt}=g_E(\textbf{x},\textbf{u},\Gamma)+g_I(\textbf{x},\Gamma), 
\end{cases}
 \end{equation}
where $g_E$ and $g_I$ represents respectively the convective and the diffusion component of the equation for the surfactant,
\begin{flalign*}
&g_E=-\nabla_\gamma \cdot (\Gamma \textbf{u}_\gamma) -2H(\textbf{x})\Gamma (\textbf{u}\cdot \textbf{n})&, \\ &g_I=\frac{1}{Pe}\nabla_\gamma^2\Gamma.&
\end{flalign*}
These two terms need to be treated differently: an implicit scheme is a good choice for the diffusion term and an explicit method for the convective term \cite{IMEX}. The system is evolved in time using the combination of the Midpoint Rule for the evolution of the drop with an Implicit-Explicit (IMEX) second-order Runge-Kutta scheme for the evolution of the surfactant concentration. To make the implicit part of the solver efficient also
for large diffusion coefficients, a preconditioner is designed taking advantages of spherical harmonics eigenfunction properties \cite{sorgentone_tornberg}. The overall scheme is adaptive with respect to drop deformation and surfactant concentration; in \cite{pallson} we showed the efficiency in minimizing the number of BIM computations while still meeting the set tolerance.\\
While simulating surfactant-covered drop in electric field, significant distortions of the point distributions representing the drops surfaces may arise, especially for long time simulations or in the case of strong electric field. This can introduce unresolvable high-frequency spherical harmonic components and lead to aliasing errors and numerical instabilities. For this reason a reparametrization procedure is absolutely needed. In order for this procedure not to ruin the overall method accuracy, we developed a spectrally accurate algorithm to ensure good quality of the surface representation also for long-time simulations \cite{sorgentone_tornberg}. It is based on an optimization procedure previously developed in \cite{veerapaneni2011} combined with the global representation of both the drop position and the surfactant concentration. We will apply it to the current numerical method at the end of every time-step; this technique is fundamental not only in the presence of high distortions due to strong electric fields, but e.g. also in the simple case of a weak field where the drop reaches a steady state but there is a non-zero tangential velocity. In these cases the reparametrization procedure avoids numerical instabilities, and it is also fundamental to ensure spectral convergence.\\
The overall procedure to simulate surfactant-covered drops in electric field is summarized in Algorithm \ref{alg:pseudo_code}.

\begin{algorithm}[h!]
  \caption{Evolution of surfactant-covered drops in electric field}
  \vspace{0.3cm}
  \begin{algorithmic}
\While{$t<T_{max}$}\\
  \vspace{0.3cm}
\textbf{Given $\textbf{x}^{(t)},\Gamma^{(t)}$:}
    \State - Compute the interfacial force $\mbf$, eq. (\ref{eq:interfacial_force})
    \State - Compute the electric force $\mbf^E$, eq. (\ref{eq:el_force}) \Comment{see steps (i)-(iv) in Section \ref{sc:mathematical_formulation}}
    \State - Compute the velocity $\textbf{u}^{(t)}$ by solving eq. (\ref{eq:BIM_galerkin})
\Large
$$\Downarrow$$
\normalsize
\begin{equation*} 
\begin{split}
\text{Compute intermediate values: } \begin{cases} 
\textbf{x}^{(t+dt/2)}=\textbf{x}^{(t)}+\frac{dt}{2}\textbf{u}^{(t)}\\ 
\Gamma^{(t+dt/2)}=\Gamma^t+\frac{dt}{2}[g_E(\textbf{x}^{(t)},\textbf{u}^{(t)},\Gamma^{(t)})+\frac{1}{Pe}\nabla_S\Gamma^{(t+dt/2)})] \\ 
\end{cases}
\end{split}
\end{equation*}
\State - Compute the interfacial force, eq. (\ref{eq:interfacial_force})
\State - Compute the electric force, eq. (\ref{eq:el_force}) 
\Comment{see steps (i)-(iv) in Section \ref{sc:mathematical_formulation}}
\State - Compute the velocity $\textbf{u}^{(t+dt/2)}$ by solving eq. (\ref{eq:BIM_galerkin})
\Large
$$\Downarrow$$
\normalsize
Based on these quantities, update:
\begin{equation*} 
\begin{split}
\begin{cases} 
\textbf{x}^{(t+dt)}=\textbf{x}^{(t+dt/2)}+dt\textbf{u}^{(t+dt/2)}\\
\Gamma^{(t+dt)}=\Gamma^t+dt[g_E(\textbf{x}^{(t+dt/2)},\textbf{u}^{(t+dt/2)},\Gamma^{(t+dt/2)})+g_I(\textbf{x}^{(t+dt/2)},\Gamma^{(t+dt/2)})]\\
\end{cases}
\end{split}
\end{equation*}
\\
\\\hrulefill 
\\
\textbf{Adaptivity in time:}
    \State $err=max(err_{drop},err_{surfactant})$
    \If {$err<tol$}
        \State Reparametrization
        \State t=t+dt;
    \EndIf \\
 $dt=dt(0.9 \frac{tol}{err})^{1/2}$
\\\hrulefill 
\\ 
 \EndWhile

   \end{algorithmic}
       \label{alg:pseudo_code}
  \end{algorithm}

\section{Small Perturbation Theory} \label{sc:spt}
In the limit of small deviations from sphericity,  drop shape is parametrized relative to the equilibrium spherical shape. The position of the interface is given by  
\begin{equation}
\mbx_s=  a \left( 1+F(\theta, \phi,t) \right) \rhat \,,\quad\mbox{where}\quad F(\theta, \phi,t)=\sum_{n=0}^{\infty}\sum_{m=-j}^j F_{nm} (t)\tilde{Y}_{nm}(\theta, \phi) 
\end{equation}
where $\tilde{Y}_{nm}=e^{\im m\phi}P_n^m \left(\cos\theta\right)$ are the unnormalized spherical harmonics defined in eq. (\ref{eq:normalized_spharm}) and $\rhat$ denotes the unit radial vector. Analytical results for the shape parameters $F_{nm}$ are obtained as perturbation expansions  in the small parameter $\Ca_E$:
\begin{equation}
\label{shapesF}
F_{nm}=\Ca_E F_{nm}^{(1)}+\Ca_E^2 F_{nm}^{(2)}+O\left(\Ca_E^3\right)\,,
\end{equation}
where $n\le 8 $ in the case of a uniform 
 {$\textbf{E}_\infty=(0,0,E_u)=E_u\nabla\left(r  Y_{10}\right)$},  and a linear (quadrupole), {$\textbf{E}_\infty=E_q(-x,-y,2z)=E_q  \nabla(r^2Y_{20})$}, electric fields. Drop shapes are axisymmetric (around the $z$ axis), thus only the $m=0$ shape modes are nonzero.  Note also that $Ca_E$ is defined based on the equilibrium surface tension, so the surfactant enters implicitly.

The shape of the drop is often described in terms of the deformation parameter $D$,  introduced by Taylor in his linear theory of drop deformation in an uniform electric field \cite{Taylor159}. 
\beq 
\begin{split}
D&=\frac{a_{||}-a_{\perp}}{a_{||}+a_{\perp}}\\
&=\frac{1}{4}\Ca_E\left(3 F_{20}^{(1)}+\frac{5}{4} F_{40}^{(1)}\right)\\
&+\frac{1}{4}\Ca_E^2\left[-\frac{3}{4}\left(F_{20}^{(1)}\right)^2-\frac{55}{64}\left(F_{40}^{(1)}\right)^2-\frac{19}{8} F_{20}^{(1)}F_{40}^{(1)}+3 F_{20}^{(2)}+\frac{5}{4} F_{40}^{(2)} +\frac{21}{8} F_{60}^{(2)}+\frac{93}{64} F_{80}^{(2)} \right]\,,
\end{split}
\label{eq:def_num} 
\eeq
 where $a_{||}$ and $a_{\perp}$ are the drop lengths in direction parallel and perpendicular to the applied field.  Due to the axial symmetry, at steady state  the electric stresses are balanced by Marangoni stresses. Accordingly, there is no flow and the stationary drop shape is independent of the viscosity ratio $\lambda$\footnote{This is true only for the surfactant-covered case and not for a clean drop; in the Appendix we list the coefficients also for the special case of a clean drop with $\lambda=1$ and we will extend the theory for general viscosity ratio in a forthcoming paper.}.

Table \ref{tab:SPT1} references the small perturbation theories (SPT) developed for uniform and quadrupole electric fields for clean and surfactant-laden drops. In the case of a uniform field, Ha and Yang \cite{Ha1995} studied the small deformation of a surfactant-covered drop with significant surfactant surface diffusion (small P\'eclet number, $Pe \approx Ca_E$), while Nganguia et al.  \cite{Nganguia2013} considered the case of a  drop covered with non-diffusing surfactant. 
In the case of quadrupole electric field and non-diffusing surfactant ($Pe=\infty$), the second order theory for drop deformation in a linear field is developed by us.  We have also reexamined the  clean drop problem. Our leading order result for the quadrupole field recovers the expressions reported by Deshmukh and Thaokar  \cite{Deshmukh_2013}, however, we found discrepancies with the second order coefficients. This issue will be further investigated in a forthcoming paper.
\subsection{Linear theory}
At leading order, the drop shape in a uniform electric field is an oblate or prolate spheroid described by \cite{Nganguia2013}
\begin{equation}
F_{20}^{(1)}=E_u^2\frac{3((1+\Rr)^2 -4Q)}{4(2+\Rr)^2},
\end{equation}
where $R$ and $Q$ have been defined before eq. (\ref{eq:EneEni}).
The electric stresses responsible for the drop deformation are quadratic in the electric field, hence at leading order drop shape in a uniform electric field 
is an ellipsoid parametrized by a single shape parameter corresponding to a harmonic of order $n=2$. However, the shape in a quadrupole electric field 
 involves two shape parameters, $n=2,4$
 \cite{Mandal_2016}:
\begin{equation}
\begin{split}
F_{20}^{(1)}&=E_q^2\frac{25\left(( 2 \Rr-1) (3 + 2 \Rr)-5Q\right)}{28\left(3+2\Rr\right)^2}\,,\\
F_{40}^{(1)}&=E_q^2\frac{10\left(1+\Rr+\Rr^2-3Q\right)}{7\left(3+2\Rr\right)^2}.
\end{split}
\end{equation}
At this order, deformation is independent of surfactant elasticity.

For the sake of completeness, here we also list the expressions for a clean drop with no viscosity contrast, $\lambda=1$. In the uniform field \cite{Taylor159},

\begin{equation}
F_{20}^{(1)}=E_u^2\frac{3( 2 + 3 \Rr + 2 \Rr^2 - 7 \Sr)}{8(2+\Rr)^2}
\end{equation}
and in the linear field \cite{Deshmukh_2012,Deshmukh_2013},
\begin{equation}
\begin{split}
F_{20}^{(1)}=E_q^2\frac{25 (-3 + 3 \Rr + 4 \Rr^2 - 4 \Sr) }{112 (3 + 2 \Rr)^2}\,,\quad F_{40}^{(1)}=E_q^2\frac{5 (4 + 3 \Rr + 4 \Rr^2 - 11 \Sr) }{56 (3 + 2 \Rr)^2}.
\end{split}
\end{equation}

\subsection{Higher-order theory}
The quadratic corrections to the shape in a uniform electric field are
 \cite{Ajayi499}:
\begin{equation}
\begin{split}
F_{20}^{(2)}&=F_{20}^{(1)}\frac{3 \left(79 \Rr^3+144 \Rr^2+\Rr (51-396 \Sr)+216 Q-94\right)}{140 (\Rr+2)^3}\,,\\
F_{40}^{(2)}&=F_{20}^{(1)}\frac{3 \left(13 R^2+28 R-54 Q+13\right)}{70 (R+2)^2}.
\end{split}
\end{equation}
As mentioned above, a useful parameter to characterize the shape of the drop is the deformation number $D$ defined in eq. (\ref{eq:def_num})
The shape parameter $D$ can be easily calculated once the shape parameters $F_{nm}$ are known and for the uniform electric field it is given by \cite{Nganguia2013}:
\begin{equation}
\label{DS2}
\begin{split}
D=\left(\frac{9((1+\Rr)^2 -4Q)}{16(2+\Rr)^2}\right)\left(\Ca_E+\Ca_E^2\left(\frac{\Rr (\Rr (139 \Rr+264)-696 Q+111)+336 Q-154}{80 (\Rr+2)^3}\right)\right).
\end{split}
\end{equation}
For the clean drop with viscosity ratio $\lambda=1$ \cite{Taylor159, Ajayi499}:
\begin{equation}
\label{DS2clean}
\begin{split}
D=\frac{9( 2 + 3 \Rr + 2 \Rr^2 - 7 \Sr)}{32(2+\Rr)^2}\left(\Ca_E+\Ca_E^2\left(\frac{3 (-308 + 293 \Rr^2 + 278 R^3 + 184\Rr - 1157 \Rr \Sr + 710 \Sr)}{640 (\Rr+2)^3}\right)\right).
\end{split}
\end{equation}
Using a similar technique as in \cite{Nganguia2013,Vlahovska:2005,vlahovska_blawzdziewicz_loewenberg_2009} we can derive\footnote{The detailed derivation for the second-order theory in quadrupole  electric field as well as the transient drop deformation is shown in a forthcoming paper.} the quadratic corrections for the quadrupole electric field. In this case the coefficients are rather lengthy so we list them in \ref{appendix1}. 
\begin{center}
\begin{table}
\begin{tabular}{T|p{1cm}|p{3cm}|p{5.2cm}}
\multicolumn{2}{c|}{ } & \textbf{Clean Drop} & \textbf{Surfactant-covered Drop }\\ 
 \hline
 & \begin{sideways} \textit{$1^{st}$ order  } \ \end{sideways} &\vspace{-1.5cm}  Taylor \cite{Taylor159} & \vspace{-1.5cm} Ha and Yang ($Pe \approx Ca_E$) \cite{Ha1995} \newline Nganguia et al. ($Pe=\infty$) \cite{Nganguia2013} \\ 
 \cline{2-4}
\multirow{2}{4em}{\textbf{UNIFORM \\}} & \begin{sideways} \textit{$2^{nd}$ order  } \ \end{sideways} &\vspace{-1.5cm} Ajayi \citep{Ajayi499} & \vspace{-1.5cm} Ha and Yang ($Pe \approx Ca_E$) \cite{Ha1995} \newline Nganguia et al. ($Pe=\infty$) \cite{Nganguia2013} \vspace{-1.5cm}\\   
 \hline
  & \begin{sideways} \textit{$1^{st}$ order  } \ \end{sideways} & \vspace{-1.5cm} Deshmukh and Thaokar \cite{Deshmukh_2012,Deshmukh_2013} & \vspace{-1.5cm} Mandal et al. \cite{Mandal_2016}\\ 
 \cline{2-4}
\multirow{2}{4em}{\textbf{QUADRUPOLE \\}} & \begin{sideways} \textit{$2^{nd}$ order  } \ \end{sideways} & \vspace{-1.5cm} Present work & \vspace{-1.5cm} Present work ($Pe=\infty$)\\   
 \hline
 \label{tab:SPT1}
\end{tabular}
\caption{Summary of small perturbation theories for a single drop in a uniform or quadrupole electric field. The surface diffusivity of the surfactant is set by the P\'eclet number ($Pe$), where $Pe=\infty$ means no diffusion.}
\end{table}
\end{center}

\section{Numerical experiments} \label{sc:num_exp}
In order to validate our code, we will compare our method with theoretical, experimental and numerical results available in literature and perform a convergence test to show the spectral accuracy of the method.
\subsection{Single drop in a uniform electric field}
 \label{sc:single_un}
We will start by validating our method in the presence of a weak uniform electric field, considering an initially spherical clean drop of radius $a=1$. 
\begin{figure}[h!]
  \centering
  \begin{tabular}[b]{c}
    \includegraphics[width=.42\linewidth]{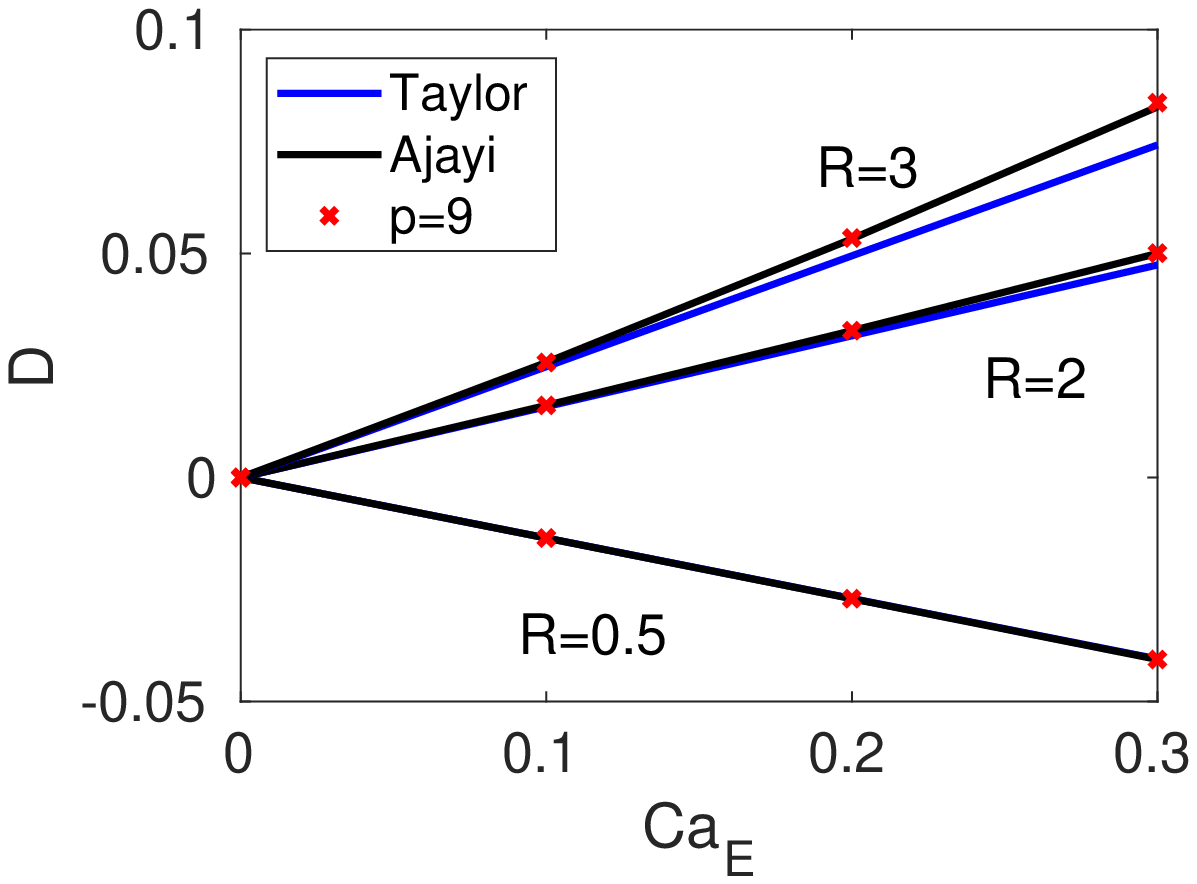} \\
    \small (a)
  \end{tabular} \qquad
  \begin{tabular}[b]{c}
    \includegraphics[width=.42\linewidth]{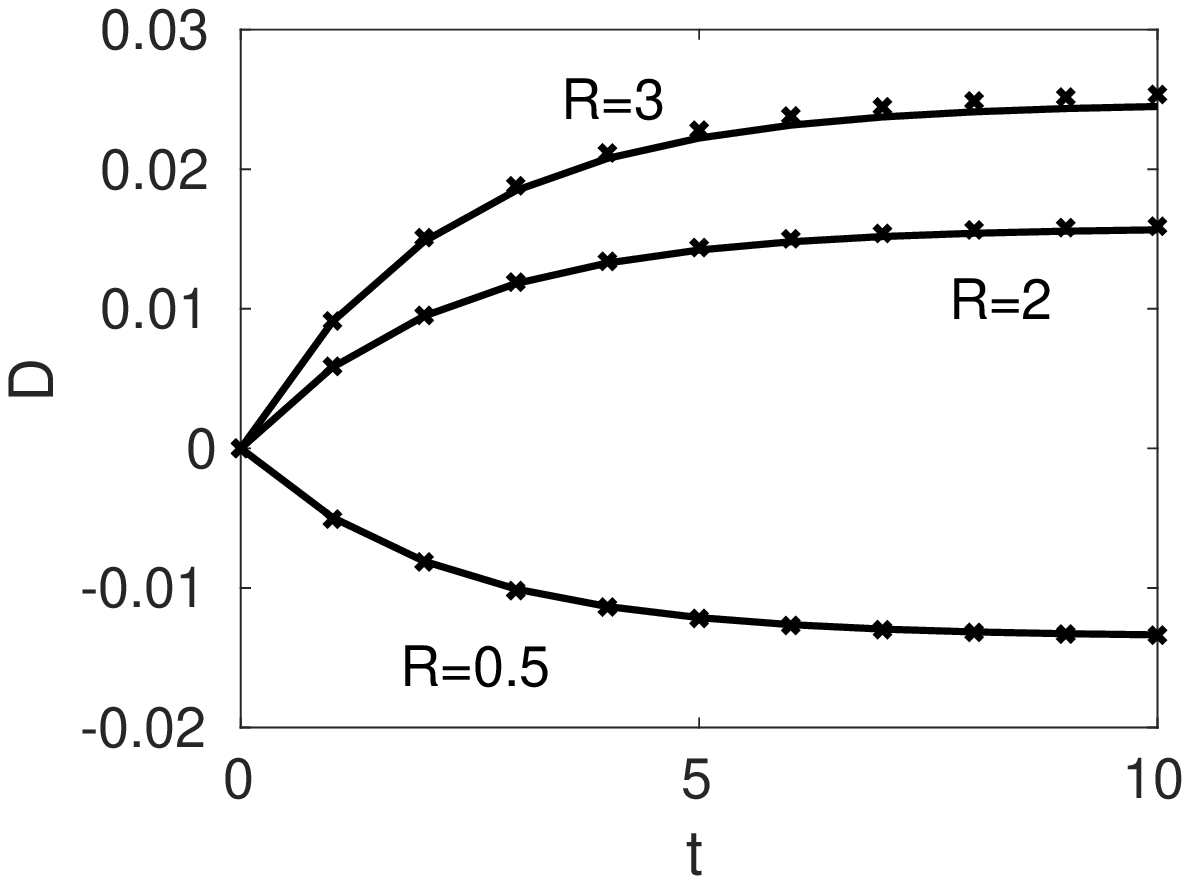} \\
    \small (b)
  \end{tabular}
  \caption{Small deformation theory for a clean drop in a uniform electric field with $Q=1$ and $\lambda=1$. a) Comparing numerical experiment (red crosses) with $1^{st}$ and $2^{nd}$ order by Taylor \cite{Taylor159} and Ajayi \cite{Ajayi499} for different values of $R$ and $Ca_E$. b) Comparing the evolution in time - eq (\ref{eq:ev_time_small_pert}) - for the case $Ca_E=0.1$, where the crosses denote the numerical simulation with $p=9$.}
	\label{fig:small_def_time}
\end{figure}

\begin{figure}[h!]
  \centering
  \begin{tabular}[b]{c}
    \includegraphics[width=.42\linewidth]{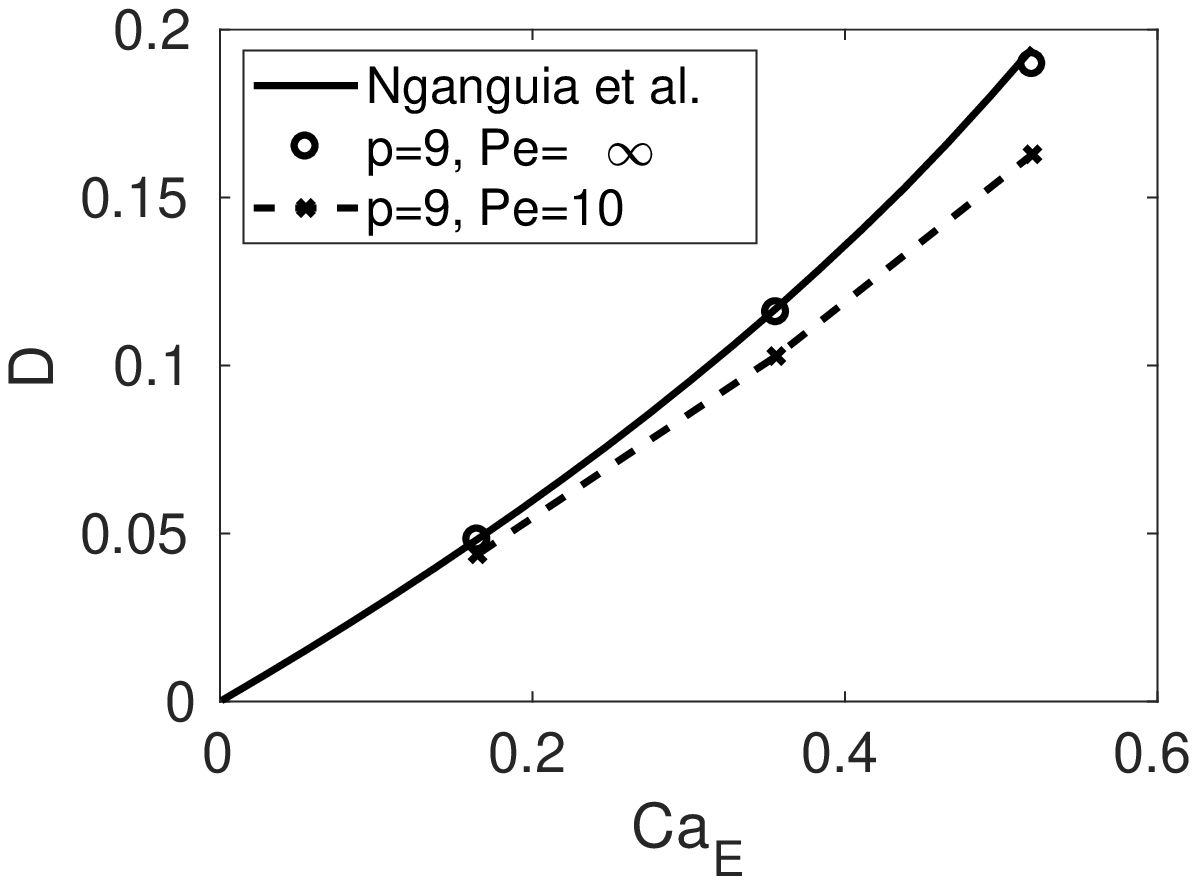} \\
    \small (a)
  \end{tabular} \qquad
  \begin{tabular}[b]{c}
    \includegraphics[width=.42\linewidth]{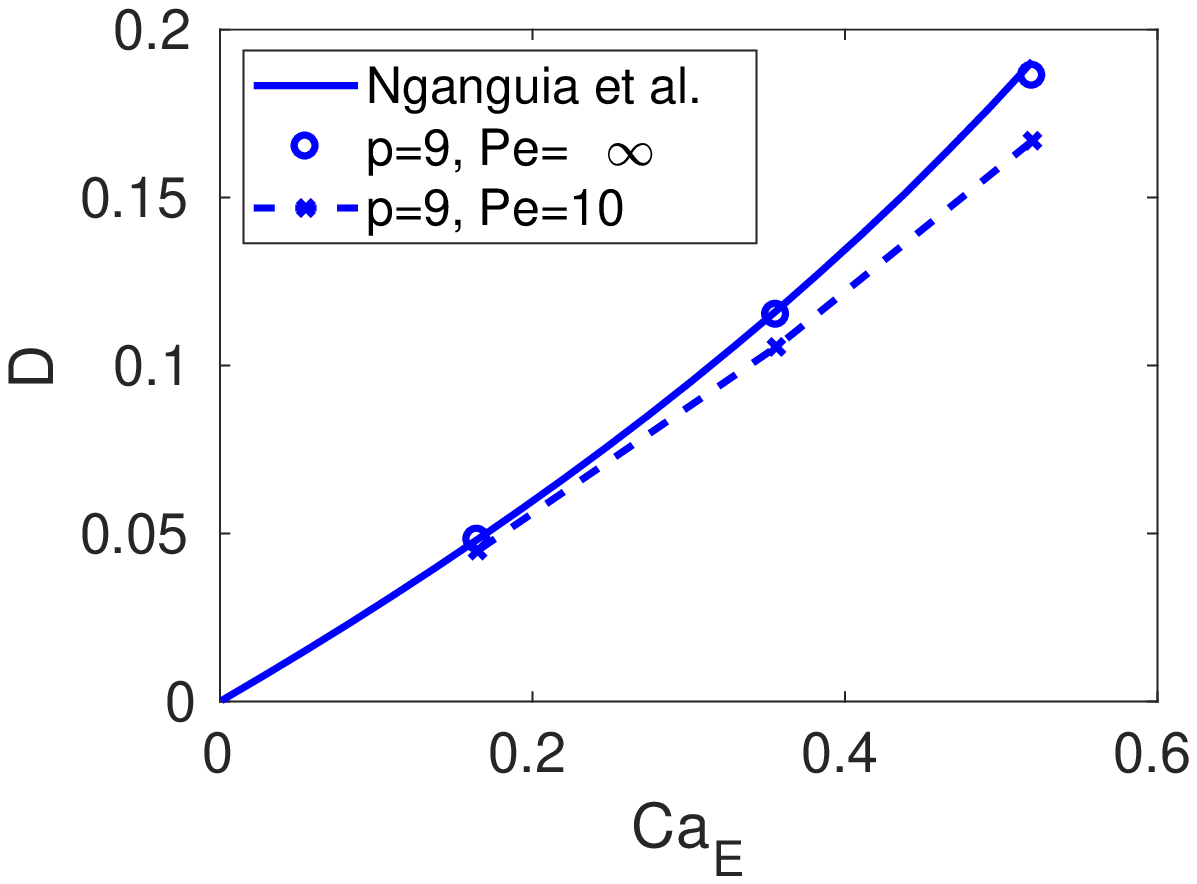} \\
    \small (b)
  \end{tabular} \qquad

  \begin{tabular}[b]{c}
    \includegraphics[width=.42\linewidth]{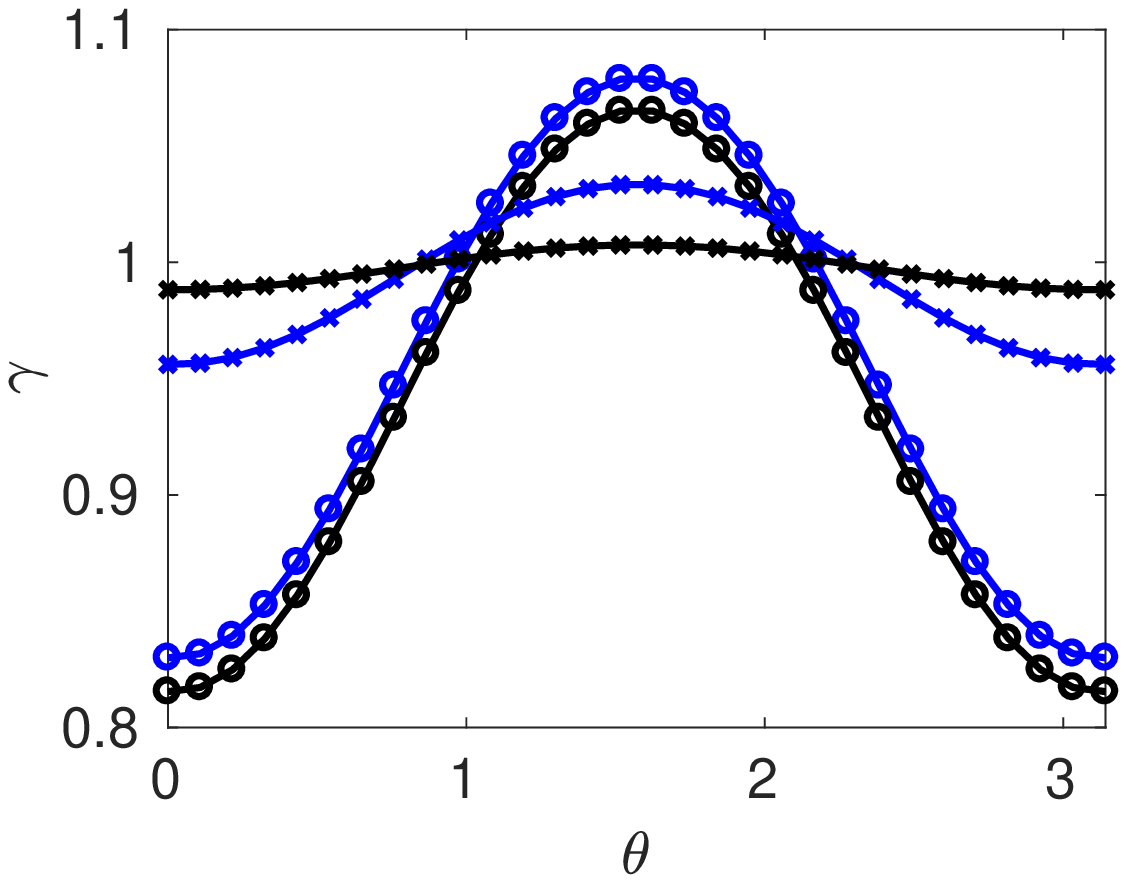} \\
    \small (c)
  \end{tabular} \qquad
  \begin{tabular}[b]{c}
    \includegraphics[width=.42\linewidth]{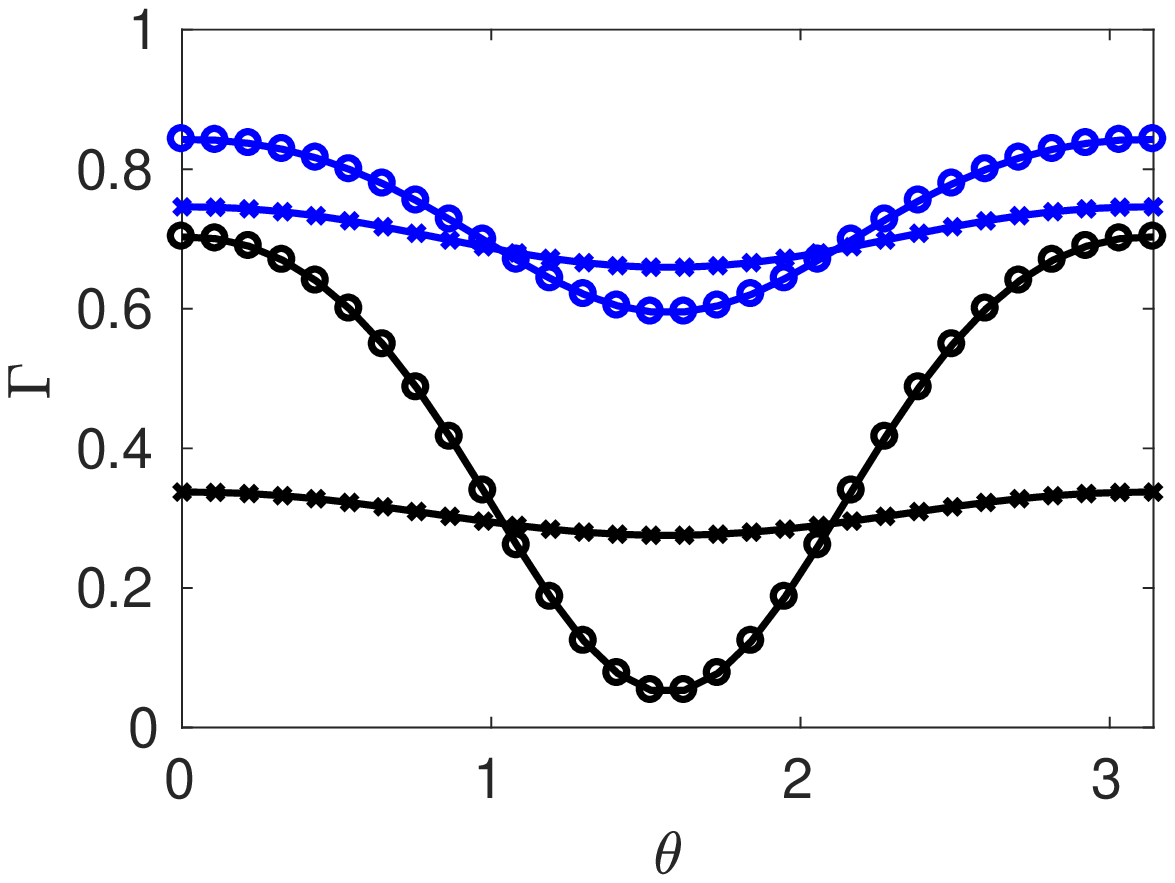} \\
    \small (d)
  \end{tabular}
  \caption{A surfactant-covered drop in a uniform electric field ($E_u=1$, $E_q=0$) with $(Q,R,\lambda)=(1,3,1)$. The deformation number is shown for a) $x_0=0.3$ and b) $x_0=0.7$, for $Pe=\infty$ (comparison with theorical results by Nganguia et al \cite{Nganguia2013}) and for $Pe=10$. The surface tension and surfactant concentration are shown respectively in c) and d) where the crosses denote the case $Pe=10$ and the circles the case $Pe=\infty$ for $x_s=0.3$ (black) and $x_s=0.7$ (blue).}
  \label{fig:D_Nganguia}
\end{figure}

As shown in Fig. \ref{fig:small_def_time}a, our results are in good agreement with the small perturbation theory, better with the second order theory then the first order one.\\
A simplified model to predict the evolution in time for small-perturbation is given by \cite{das_saintillan_2017bis}:
\beq
D(t)=D_T(1-e^{-t/\tau_D})
\label{eq:ev_time_small_pert}
\eeq
where $\tau_D= \frac{a_r\nu}{\gamma}\frac{(19\lambda+16)(2\lambda+3)}{40(\lambda+1)}$ and $D_T$ is the steady first-order deformation parameter obtained by Taylor \cite{Taylor159}. It is a good test to check the time-scaling of our numerical method, which is in good agreement with the theory as shown in Fig. \ref{fig:small_def_time}b. \\
We will now proceed to test the method on a surfactant-covered drop and compare our results with the spheroidal model by Nganguia et al. \cite{Nganguia2013} valid also for higher deformation. Note that in \cite{Nganguia2013} the diffusion is neglected ($Pe=\infty$); for completeness we also run the experiment in Fig. \ref{fig:D_Nganguia} with $Pe=10$ to compare it with the numerical results by Teigen and Munkejord, Fig. 6 in \cite{Teigen2010}. The conductivity and permittivity ratio are set as $R=3$ and $Q=1$ with $\lambda=1$ and the spherical harmonics order for the numerical experiment is set to $p=9$, which shows good agreement with both the theoretical and numerical results mentioned.

\begin{figure}[h!]
	\centering
			\includegraphics[width=75mm]{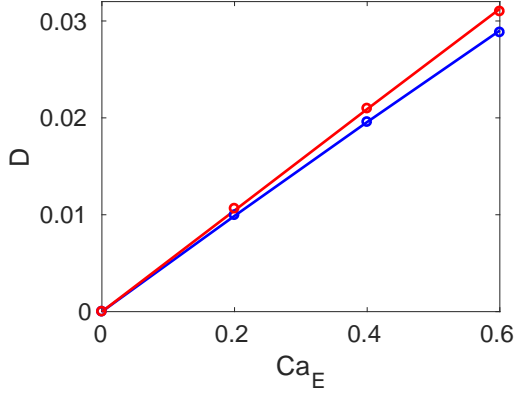}
	\caption{Comparing numerical experiment using $p=9$ (circles) with SDT by Ha and Yang \cite{Ha1995} (solid lines) for a surfactant-covered drop in a uniform electric field ($E_u=1$, $E_q=0$). The deformation number is shown for for $\lambda=1$ (blue) and $\lambda=0.01$ (red). The other parameters are $R=3$, $Q=3.5$ for the electric field and $Pe=0.4$ and $\beta=0.2$ for the surfactant.}
	\label{fig:ha_yang}
\end{figure}

The highly-diffusive case is shown in Fig. \ref{fig:ha_yang} for two different values of the viscosity ratio, $\lambda=1$ and $\lambda=0.01$. Here we compare the numerical experiment with the small perturbation theory developed by Ha and Yang \cite{Ha1995} that follows very close our numerical result.

\begin{figure}[h!]
  \centering
  \begin{tabular}[b]{c}
    \includegraphics[width=.42\linewidth]{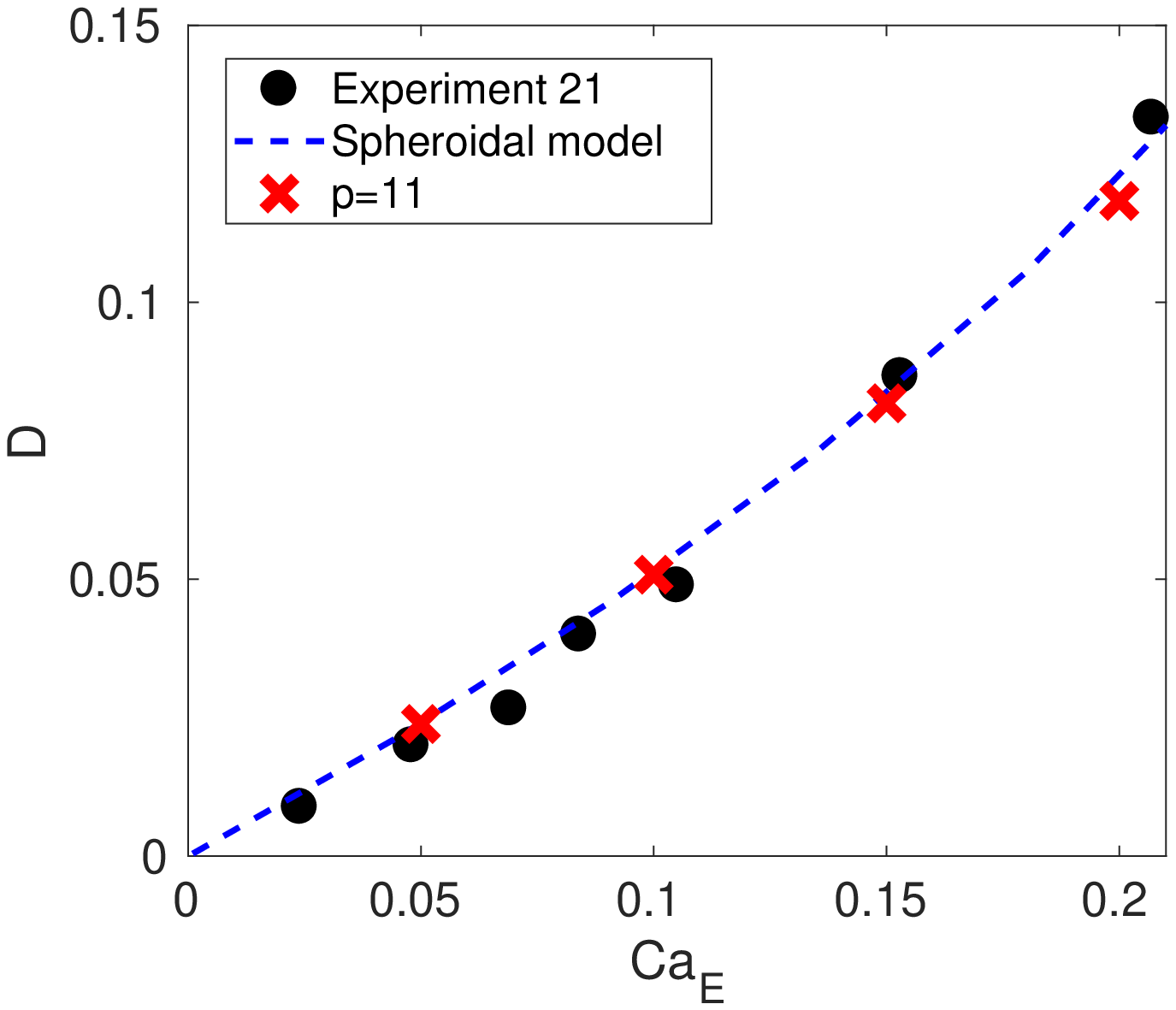} \\
    \small (a)
  \end{tabular} \qquad
  \begin{tabular}[b]{c}
    \includegraphics[width=.42\linewidth]{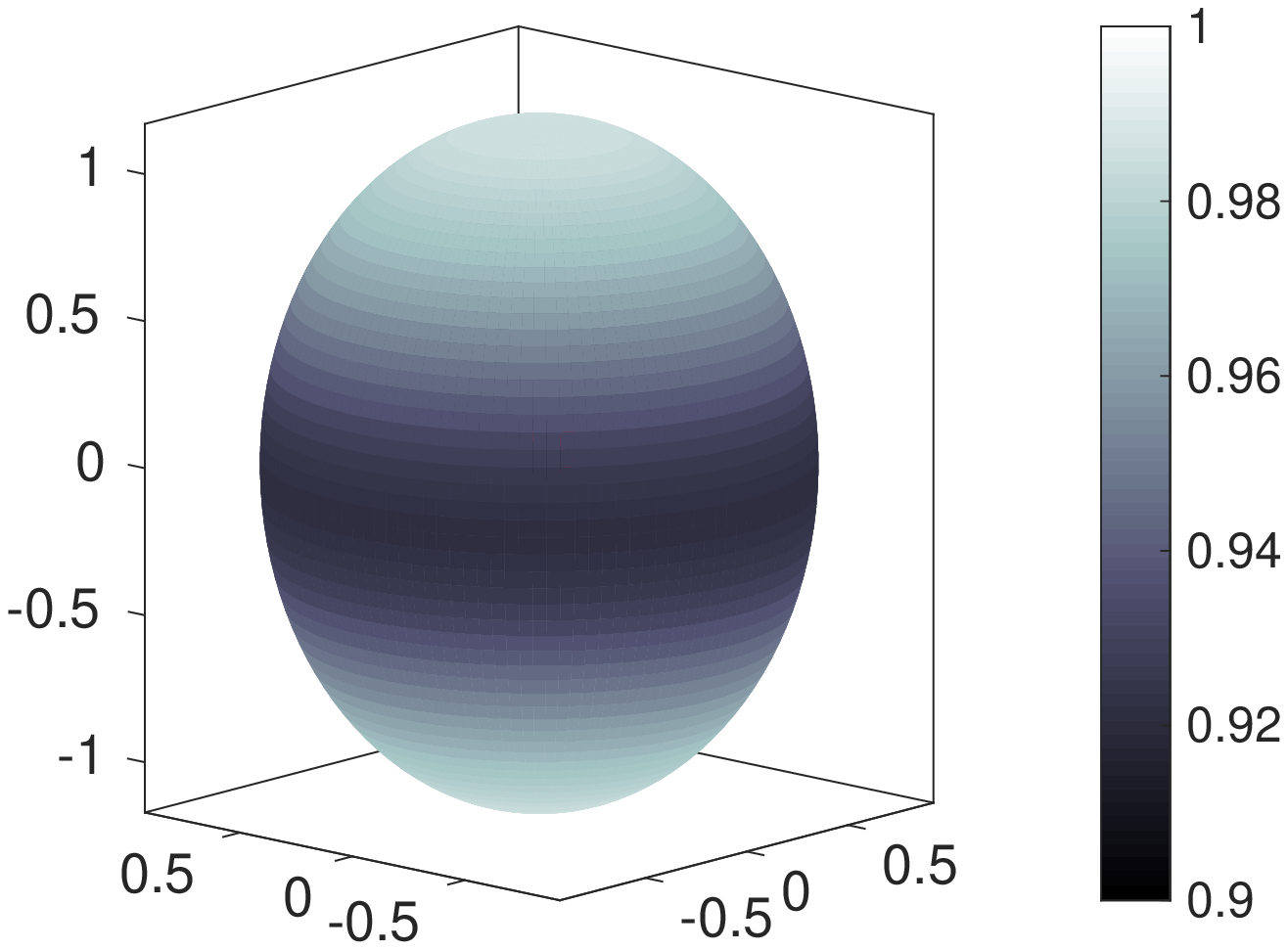} \\
    \small (b)
  \end{tabular}
  \caption{a) Comparison with spheroidal model \cite{Nganguia2013} and experimental results \cite{HA1998195} for a non-diffusive surfactant-covered drop in a uniform electric field ($E_u=1$, $E_q=0$) with $(Q,R, \lambda) = (1.3699, 10, 0.08)$ and $(\beta, x_s) = (0.04,0.96)$; b) The surfactant-covered drop at steady-state for $Ca_E=0.2$, the greyscale denotes the surfactant concentration.}
	\label{fig:nganguia3c}
\end{figure}

Ha and Yang performed also laboratory experiments with surfactant-covered drops, with a high surfactant coverage value ($x_s=0.96$) and a low elasticity number ($\beta=0.04$). Our method shows good agreement also in this case as shown in Fig. \ref{fig:nganguia3c}a; the deformed surfactant-covered drop is shown in Fig. \ref{fig:nganguia3c}b, while in Fig. \ref{fig:nganguia3c}a we show the comparison between our numerical experiments with both the experimental data and the spheroidal model previously mentioned.

\subsection{Single drop in a quadrupole electric field}
\label{sc:comp_spt_nonuniform}
\begin{figure}[h!]
  \centering
  \begin{tabular}[b]{c}
    \includegraphics[width=.42\linewidth]{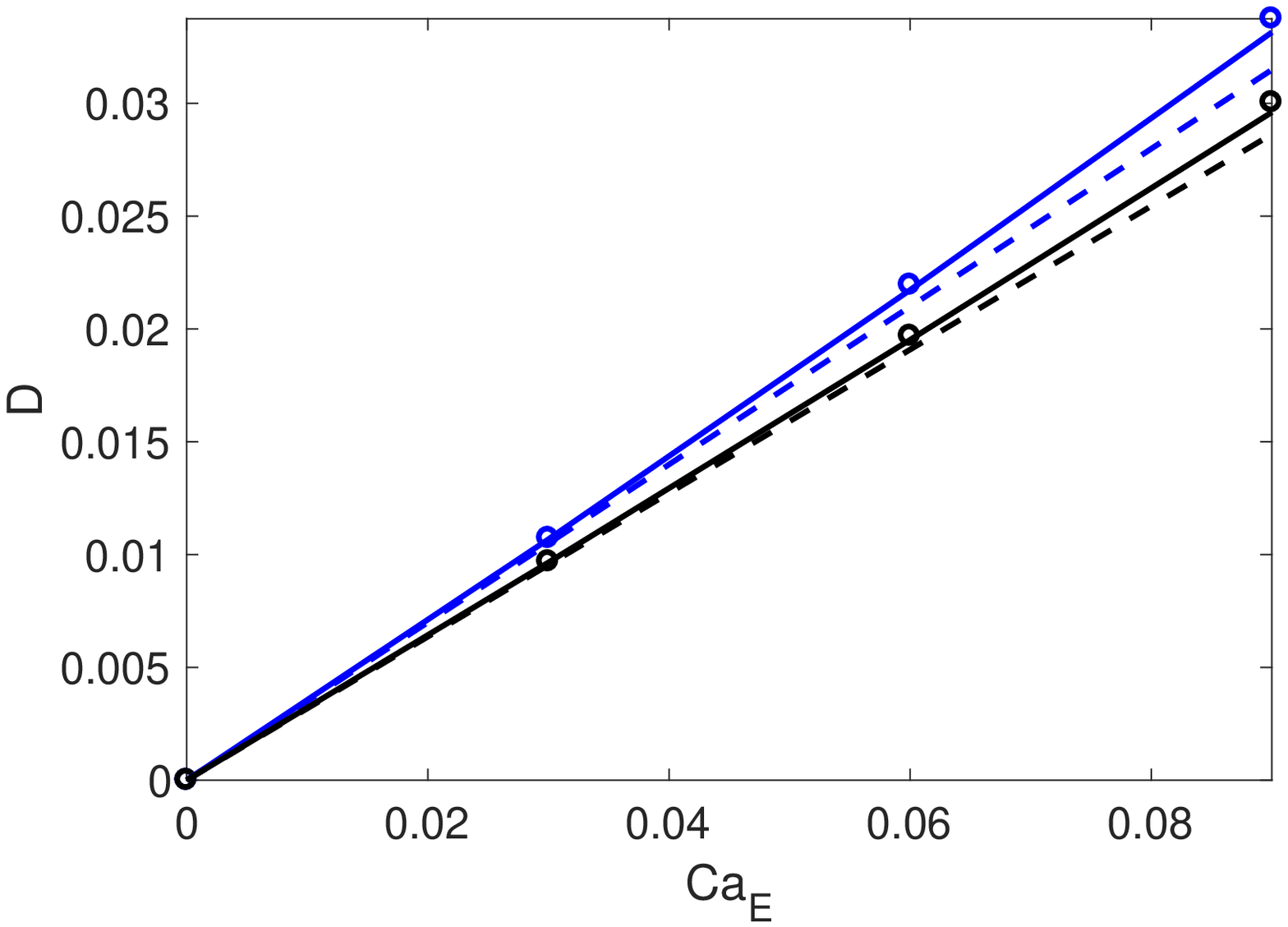} \\
    \small (a)
  \end{tabular} \qquad
  \begin{tabular}[b]{c}
    \includegraphics[width=.42\linewidth]{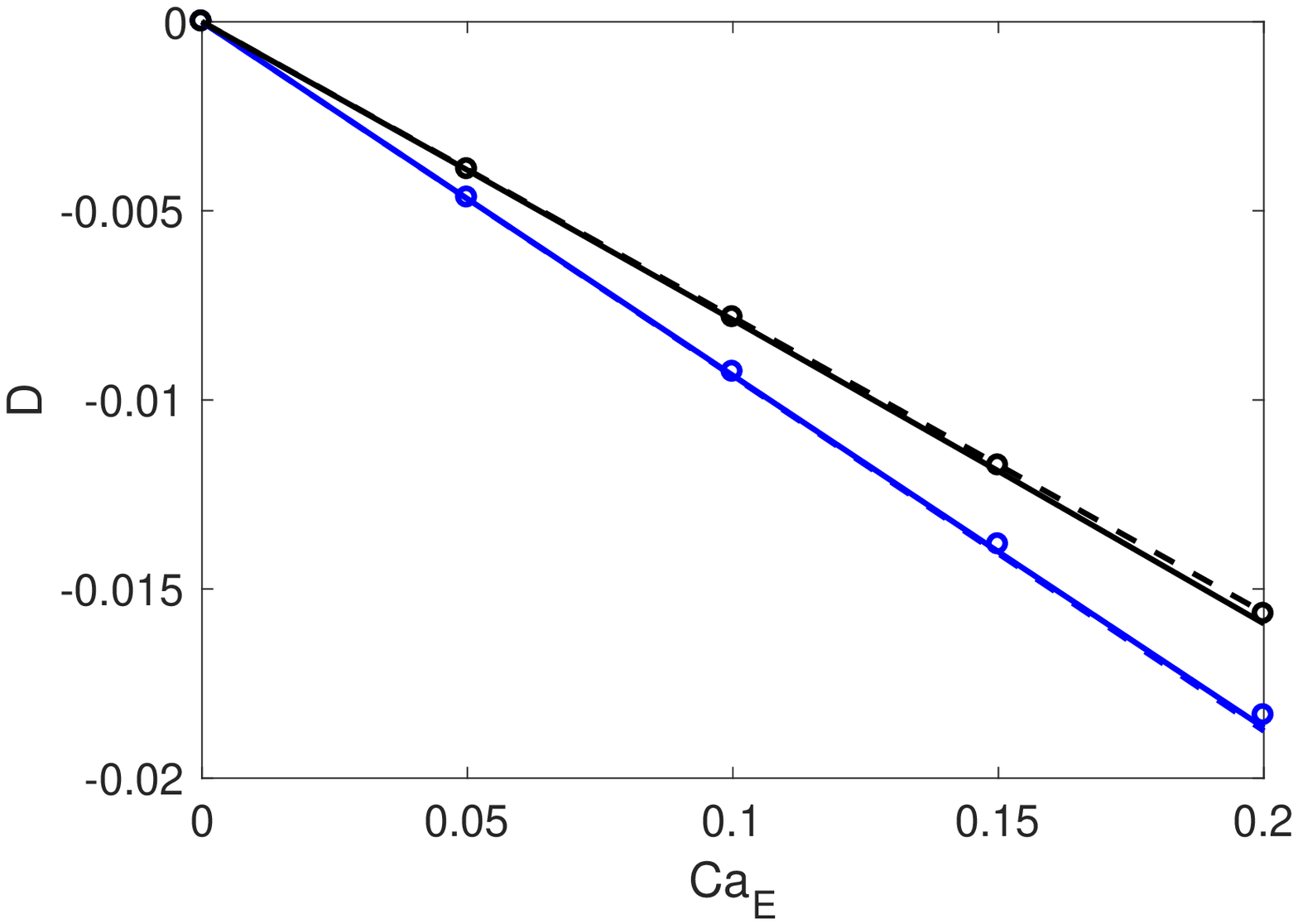} \\
    \small (b)
  \end{tabular}
  \caption{Comparison with first (dotted lines) and second-order (solid lines) SDT for a surfactant-covered drop in a quadrupole electric field presented in section \ref{sc:spt}, $E_u=0$, $E_q=1$. The circles represent the numerical simulation obtained with $p=13$ for the clean drop (black) and the surfactant-covered drop (blue) with no diffusion. a) $R=2$, $Q=0.01$, $\beta=1$ b) $R=1$, $Q=1.5$, $\beta=0.5$.}
	\label{fig:petia-2nd}
\end{figure}

We will now test our method against the second order small deformation theory developed in section \ref{sc:spt} for surfactant-covered drops in a quadrupole electric field. We ran several test cases and once again we find good agreement between the theory and the numerical simulations. Fig. \ref{fig:petia-2nd}) show two different simulations comparing the deformation of the clean and surfactant covered drop. To better visualize how the first and second order theories differ, we also show in Tab. \ref{table:quad_sdt_diff} the relative difference in the deformation number obtained by numerical simulations as compared to these theoretical results for the case shown in Fig. \ref{fig:petia-2nd}a. As expected, the best agreement is obtained for the smallest capillary number, as the small perturbation theory gets less accurate as the capillary number is increased. That said, the second order theory clearly matches the numerical results better than the first order theory for all the capillary numbers considered.\\ For the physical parameters set in Fig. \ref{fig:petia-2nd}b the second order contribution is not so important, but a sharper difference between the clean and the surfactant-covered case is shown. \\
Observe also that the surfactant-covered drop cases have been run for three different values, $\lambda=0.1,1,2$ but the deformation number is not affected by the viscosity ratio as already discussed in section \ref{sc:spt}, so the three simulations coincide; the only difference is the time needed for reaching the steady state that depends on $\lambda$. For these cases we also checked numerically that the tangential component of the velocity is zero when the steady state is reached, which is not true for the clean case, where the tangential velocity is kept active even when the steady state is reached.

\begin{table}[h]
\begin{center}
 \begin{tabular}{l | C | C  | C |}
 \hline 
 \multicolumn{1}{|c||}{$Ca_E$} & 0.03 & 0.06 & 0.09   \\
            \hline \hline \hline
             & \multicolumn{3}{|c|}{\textbf{Clean drop}} \\
            
            \hline
            \multicolumn{1}{|c|}{1st order}    & 1.4696\% &   3.0262\% &   4.6807\% \\
            \multicolumn{1}{|c|}{2nd order}   & 0.3517\%  &  0.8255\%  &  1.4361\%  \\
          \hline  \hline 
          
\cline{2-4}
             & \multicolumn{3}{|c|}{\textbf{Surfactant-covered drop}}\\
           
            \hline
            \multicolumn{1}{|c|}{1st order}    & 2.1193\% &    4.3641\% &    6.7063\% \\
            \multicolumn{1}{|c|}{2nd order}   & 0.3972\% &    0.9989\% &    1.7821\%  \\  
          \hline 
           
        \end{tabular}

\caption{Relative difference between the numerical and the theoretical values for the deformation number $D$ of the experiment shown in Fig. \ref{fig:petia-2nd}a. The theoretical values from small perturbation theory (first/second order) are only expected to be accurate for small electric capillary numbers.}
 \label{table:quad_sdt_diff}
 \end{center}
\end{table}


\subsection{Convergence test for drop position and surfactant distribution}
\begin{figure}[h!]
  \centering
  \begin{tabular}[b]{c}
    \includegraphics[width=.42\linewidth]{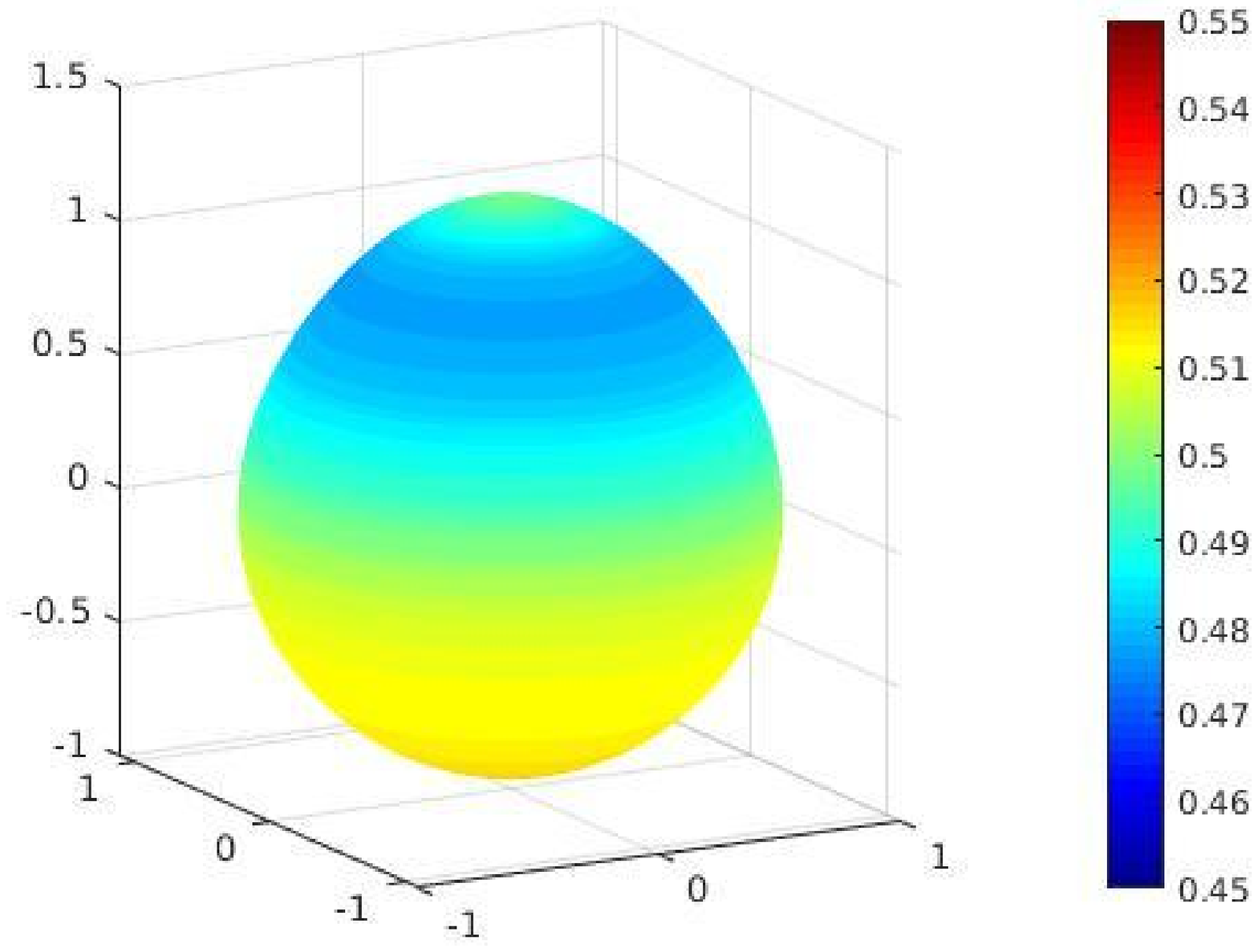} \\
    \small (a)
  \end{tabular} \qquad
  \begin{tabular}[b]{c}
    \includegraphics[width=.42\linewidth]{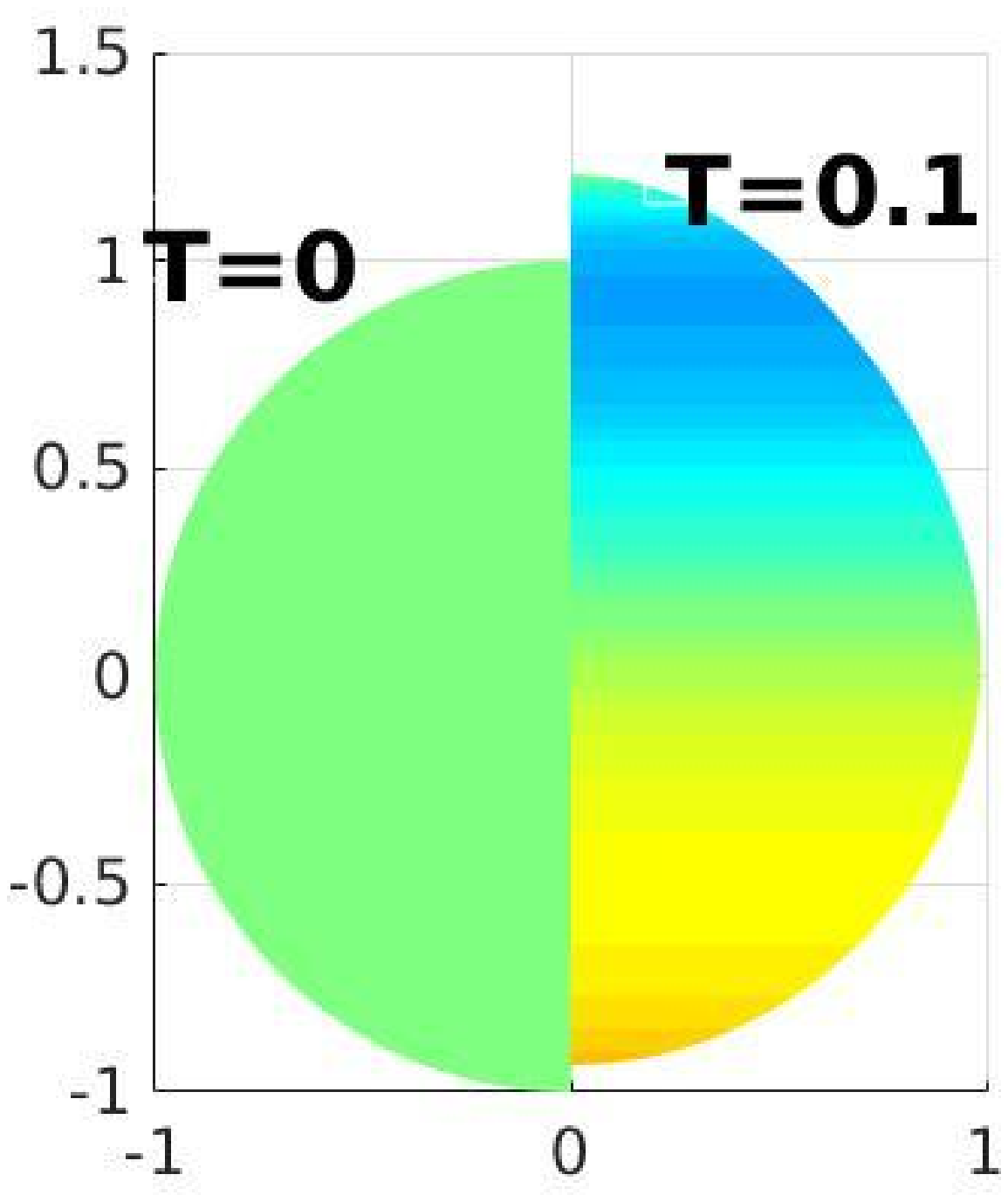} \\
    \small (b)
  \end{tabular}
  \caption{The combined effect of linear and quadrupole fields on a single surfactant-covered drop with $R=6$, $Q=2$, $Pe=100$, $\beta=0.2$, $x_s=0.5$ and $\tilde{\beta}=1$. a) The surfactant-covered drop at time T=0.1. The colorbar denotes the surfactant concentration. (b) A section of the drop, comparing the shape and the surfactant concentration at time T=0 and T=0.1.}
	\label{fig:drop_conv_comp}
\end{figure}

\begin{figure}[h!]
	\centering
			\includegraphics[width=65mm]{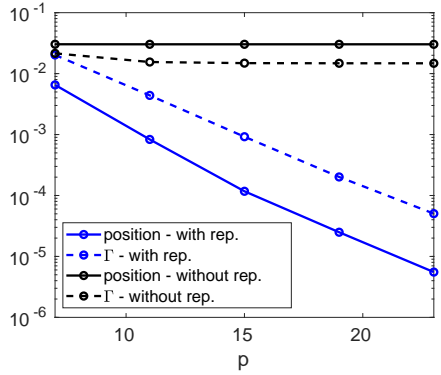}
	\caption{$L_\infty$ error for the position vector (solid line) and for the surfactant concentration (dotted line) for the test case $R=6$, $Q=2$, $Pe=100$ and $\tilde{\beta}=1$.}
	\label{fig:convergence}
\end{figure}

In Fig. \ref{fig:convergence} it is shown the convergence of our method for a test case where we consider the combination of linear and quadrupole fields with $R=6$, $Q=2$, $Pe=100$, $\beta=0.2$, $x_s=0.5$ and $\tilde{\beta}=1$. The time-step tolerance is set to $tol=10^{-7}$ and the reference solution is obtained with $p=31$. The shape and surfactant concentration of the drop at the final time are shown in Fig. \ref{fig:drop_conv_comp}a whilst in Fig. \ref{fig:drop_conv_comp}b we compare a section of the drop at initial and final times. We can observe that, in this special case, the drop is not only deforming but also translating. This is due to the combined effect of a linear and a quadrupole electric field. The translation sweeps the surfactant towards the back (bottom) of the drop, while the electrohydrodynamic flow pushes the surfactant towards the drop poles (in the considered case $R/Q>1$). This gives rise  to the peculiar behavior that we see in Fig. \ref{fig:drop_conv_comp}, with a slight increase in the surfactant concentration at the front (top) of the drop, in addition to the surfactant accumulation at the drop bottom.

In Fig. \ref{fig:convergence} we show the $L_\infty$ error for the same simulation run with and without reparametrization: it is clear that the spectral convergence is immediately lost, even for such a short numerical simulation if we remove the special reparametrization procedure described in Section \ref{sc:time_and_rep}.

\subsection{Multiple drops}
\begin{figure}[h!]
	\centering
			\includegraphics[width=75mm]{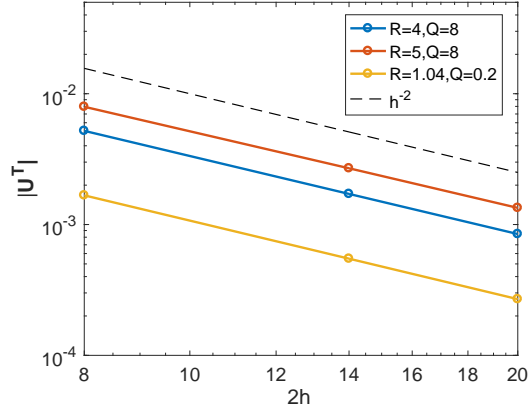}
	\caption{Translation velocity for a drop pair vs centroid separation, comparison with the predicted slope by Baygents et al. \cite{baygents_rivette_stone_1998}.}
	\label{fig:baygents-2}
\end{figure}

In this section we want to see how multiple drops interact in the presence of an electric field. Baygents et al. \cite{baygents_rivette_stone_1998} predicted how the relative velocity for a drop pair should scale with respect to the centroid separation. For the leaky-dielectric case, they predict that the translation velocity $U^T$ should scale as:
\beq
U^T=O(\frac{\tilde{\epsilon}\epsilon_0E_\infty^2a^3}{\mu h^2}),
\eeq
where $h$ is the distance between the two drops centroids. Note that the previous expression doesn't take into account how the electric parameters affect the scaling, but only how the distance between drops affect the translation velocity. We performed an experiment with two drops aligned on the $z$-axis and we checked that the scaling is correct, see Fig. \ref{fig:baygents-2}.

\begin{figure}[h!]
  \centering
    \hspace{-2cm}
  \begin{tabular}[b]{c}
    \includegraphics[width=.75\linewidth]{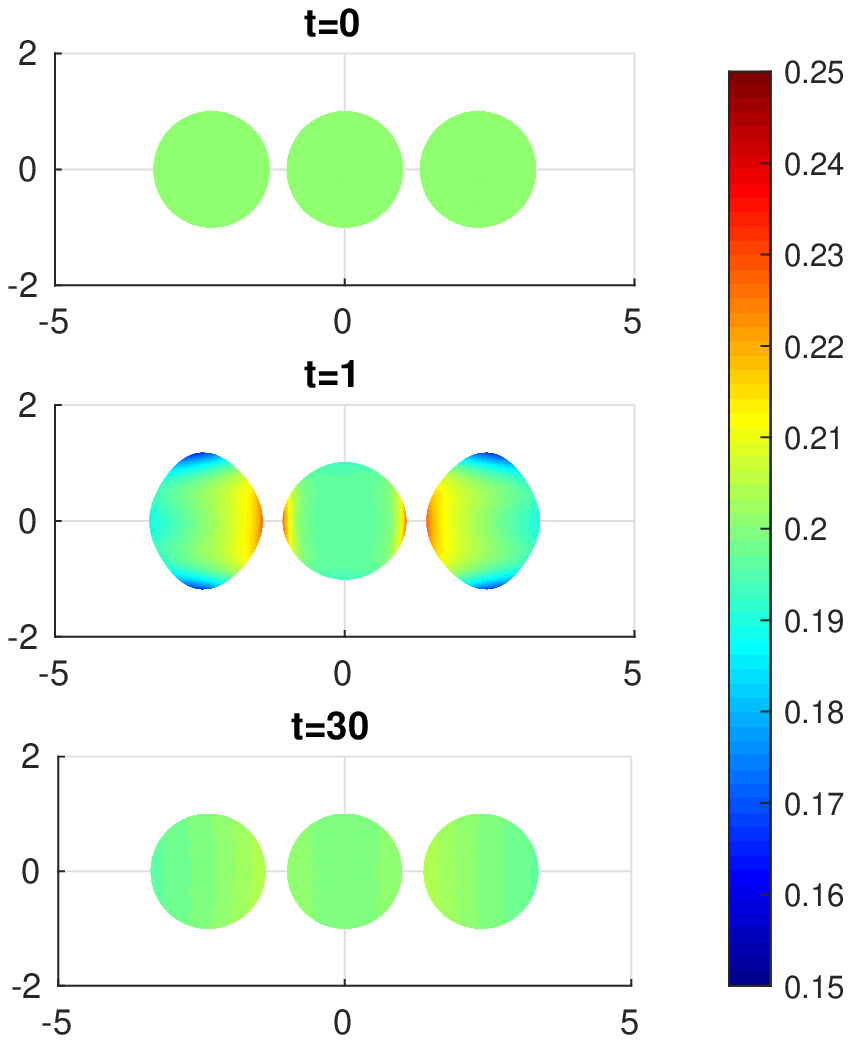} \\
    \small (a)
  \end{tabular} \qquad
  \hspace{-1.5cm}
  \begin{tabular}[b]{c}
   \includegraphics[width=.25\linewidth]{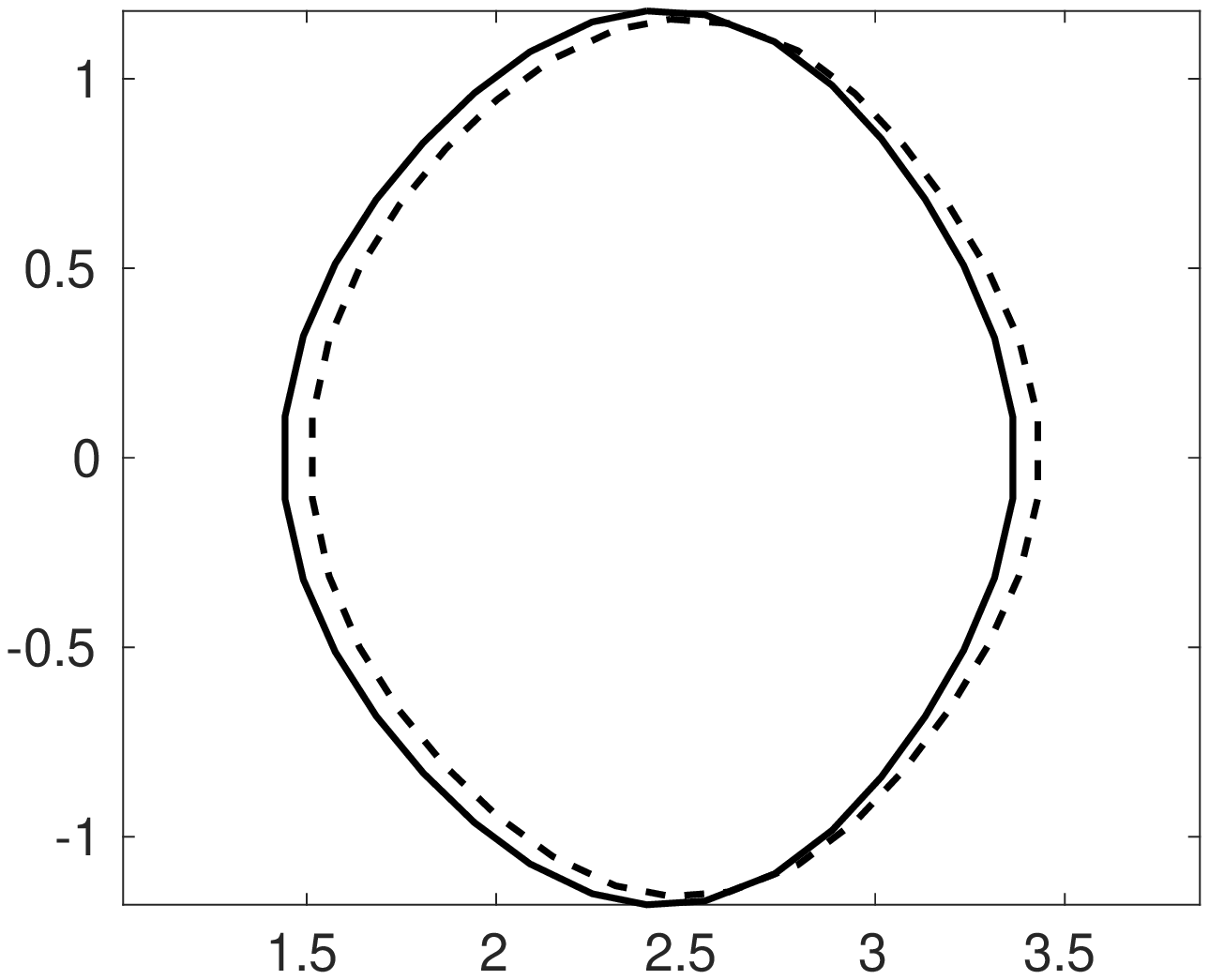} \\
    \small (b)\\
        \vspace{-0.3cm}
    \includegraphics[width=.33\linewidth]{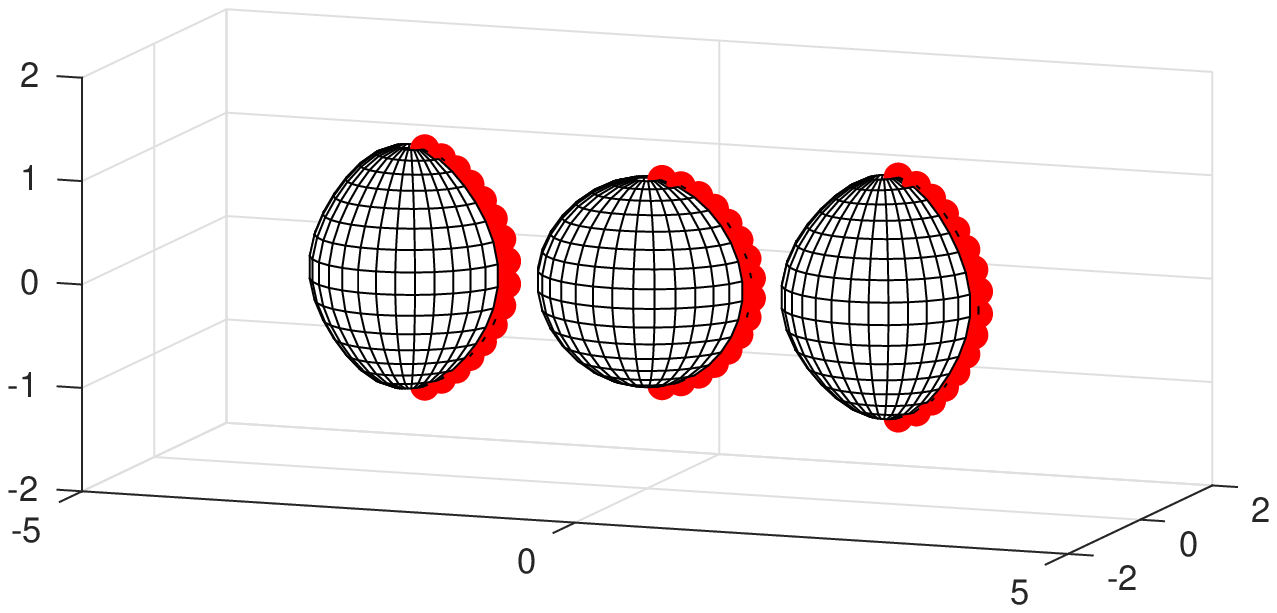} \\
    \small (c)\\
     \includegraphics[width=.25\linewidth]{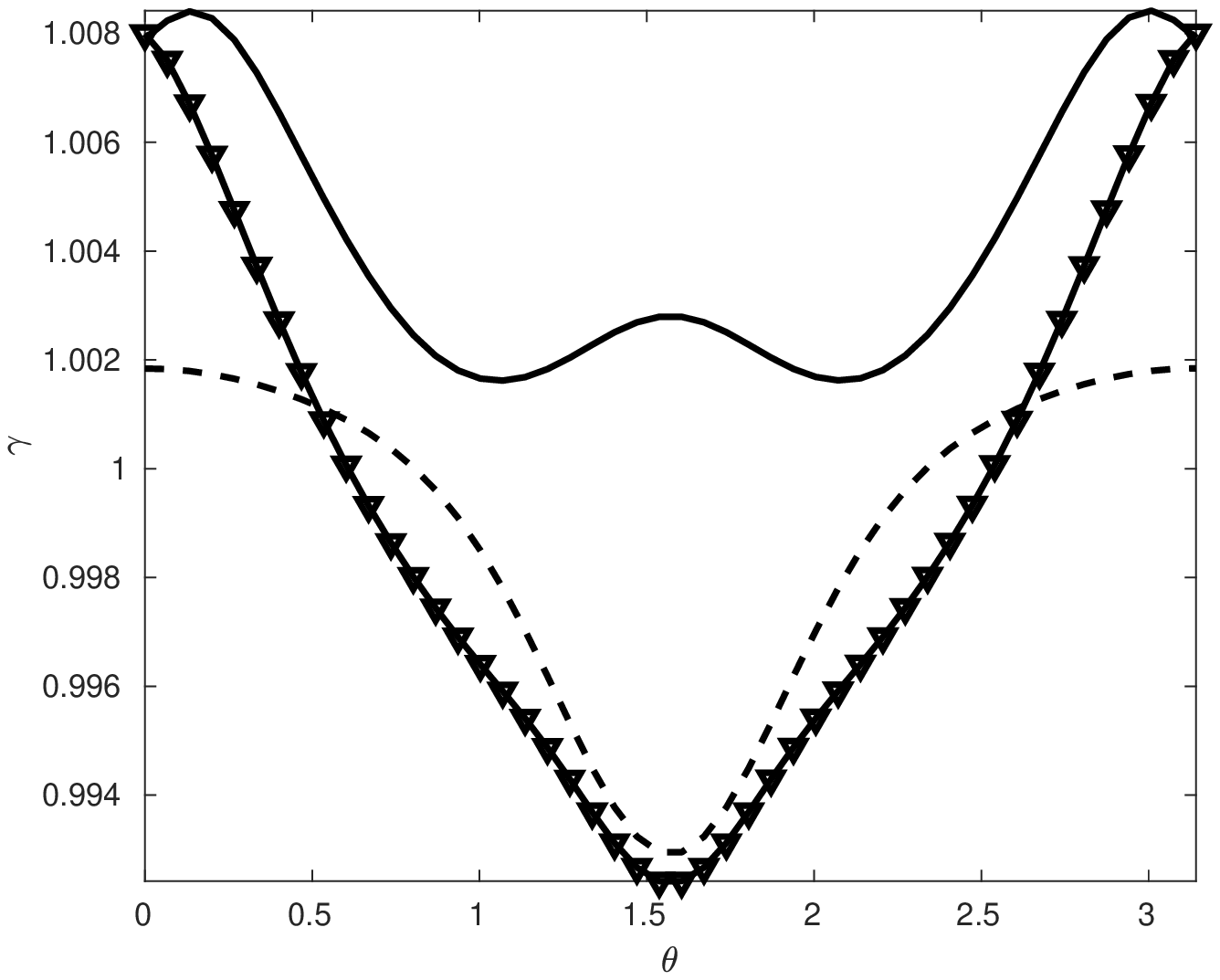} \\
    \small (d)
  \end{tabular}
  \caption{Three surfactant-covered drops interacting in a non-uniform electric field with $(Q,R,\lambda)=(5,4.5,1)$, $Ca_p=0, Ca_Q=1, Ca_E=0.1$ until time $T=1$, then $Ca_E=0$. a) The surfactant-covered drops at different times t=0,1,30. b) Splitting the effect of the quadrupole field and the interaction between drops at time $t=1$. The solid line is the contour of the drop on the right in a) but without surfactant, the dotted line is the same clean drop but for the simulation run with no other drops present. c) The points corresponding to $\phi=0$, used to show the surface tension in d), for the drop on the right (solid line), in the middle (dotted line) and on the left (v-line).}
	\label{fig:3drops}
\end{figure}

We now perform a simulation where three aligned drops on the $x$-axis are interacting in a non-uniform electric field. The main purpose of this numerical simulation is to show the robustness of our method, but of course further studies on multiple interacting drops are needed in order to better understand the electrohydrodynamics of these systems, and we will extend this kind of analysis in the near future.
For this simulation $R=4.5, Q=5$ up to time T=1, when the electric field is switched off and the particles relaxed to their spherical shape with uniform surfactant concentration. 
In Fig. \ref{fig:3drops}a) we can see how the three particles react differently to the quadrupole field. To further understand how the interaction between particles and how the non-uniform field are affecting the simulation, we consider the same experiment for clean drops (so that the surfactant is not interfering), for which we run the simulation both for the three particles and also for a single particle, the particle on the right. We compare the results in Fig. \ref{fig:3drops}b), where we can see that the shape of the drop is due to the quadrupole field that act differently depending on the drop position. The interaction of the three drops result in an attractive effect between the particles (we can see the translation of the solid line in Fig. \ref{fig:3drops}a)).\\
Setting the time-step tolerance to $tol=1e-06$ and the order of the spherical harmonics expansion to $p=15$, we obtain a relative error in surfactant conservation (measured between time zero and the final time $T=30$) lower than $3e-05$ for all the three particles, and a relative error for the volume less than $7e-08$ for both the clean and the surfactant-covered case. 

\section{Conclusions and future work}
\label{conclusions}
We investigated the behavior of surfactant-covered drops in electric field both theoretically and numerically. We introduced a new second-order small deformation theory for the deformation of a drop placed in a weak quadrupole field in the presence of insoluble surfactant. In order to study also higher deformation and multiple drop configurations, we introduced a 3D highly accurate numerical method based on a boundary integral formulation. We validated the method against theoretical, experimental and numerical results, in all cases obtaining good agreement with the existing literature. The high accuracy of the method is achieved thanks to the special numerical tools used for the quadrature and for the reparametrization needed for maintaining a good quality of the surface representations. The method is able to handle different viscosity ratios and close interactions of drops; no simplification in terms of axysymmetric configuration is assumed, so we are able to handle multiple drops in a general setting. The adaptivity in time together with the spherical harmonics based representation makes the code efficient for a single drop or for a limited number of drops; fast methods have already been developed in our group \cite{fast1,fast2,fast3} and can be implemented in case of large numbers of particles.

With the present paper we have introduced and validated the numerical method which can be now used for more extensive numerical investigations to answer different fundamental physical questions. In terms of further numerical developments we would like to include geometries (e.g. walls, solids, etc) and, as already mentioned in the introduction, consider also the equation for the surface charge conservation. Including charge convection is challenging due to the non-linearity of the coupling between the charge distribution and the resulting fluid flow as shown by Das and Saintillan \cite{das_saintillan_2017bis,das_saintillan_2017}. It is a complicated subjected to study numerically due to the unsteady chaotic dynamics observed experimentally \cite{Salipante_vlahovska} and to the charge shocks that may arise in these kind of simulations \cite{lanauze_walker_khair_2015}, but it would be very interesting to see how including the surface charge would affect the multiple drop interactions.

\appendix
\section{Second-order small deformation theory for a drop in a quadrupolar electric field: the shape parameters}
\label{appendix1}
\subsection {Surfactant-covered drop with $Pe=\infty$}
\begin{equation}
\begin{split}
F_{20}^{(2)}&=\left[{6830208 \tilde{\beta} (2 R+3)^5 (4
   R+5)}\right]^{-1}\left\{125 s^4 \left(\tilde{\beta} \left(R \left(8072507-2R \left(796000 R^2+421836
   R-3212985\right)\right) Q \right.\right.\right.\\
   &\left.\left. \left.+198 (R (3124 R-15439)-24000) Q^2 \right.\right.\right.\\
   &\left.\left. \left.+R (2 R (2 R (2 R (16 R
   (6619 R+20507)+78307)-1450283)-3261593)-469541)\right.\right.\right.\\
   &\left.\left. \left.+615 (3769 Q+2447)\right)-12 (2 R+3) (4
   R+5) (R-Q) \left(36 R^2-710 R+1001 Q-327\right)\right)\right\}\, ,\\
F_{40}^{(2)}&=\left[{3699696 \tilde{\beta} (2 R+3)^5 (4
   R+5)}\right]^{-1}\left\{5 s^4 \left(\tilde{\beta} \left(135 (16 R (17273 R+15407)-134085) Q^2\right.\right.\right.\\
  &\left.\left. \left.-R (R (4 R
   (8850472 R+20123595)+48389841)-23699929) Q\right.\right.\right.\\
   &\left.\left. \left.+R (R (8 R (R (R (998636
   R+3648867)+5130720)+1972471)-19371375)-18396837)\right.\right.\right.\\
   &\left.\left. \left.+35611530 Q-3750015\right)+30 (2 R+3)
   (4 R+5) (8 R (120 R+431)-4081 Q-327) (R-Q)\right)\right\},\\
F_{60}^{(2)}&=\left[{47811456 \tilde{\beta} (2 R+3)^4 (4 R+5)}\right]^{-1}\left\{25 s^4 \left(\tilde{\beta} \left(-R \left(26872424 R^2+53198046
   R+32804787\right) Q\right.\right.\right.\\
   &\left.\left.\left.+27 (1109908 R+916985) Q^2+R (2 R (4 R (R (700708
   R+2159279)+2778484)+6293339)+868415)\right.\right.\right.\\
   &\left.\left.\left.+5 (684607 Q-767350)\right)+546 (4 R+5) \left(84
   R^2+338 R-407 Q-15\right) (R-Q)\right)\right\}\,,\\
F_{80}^{(2)}&=\left[{16648632 \tilde{\beta} (2 R+3)^4}\right]^{-1}\left\{25 s^4 \left[\tilde{\beta}(R (98888 R+121059)-318835 Q+98888) \left(R^2+R-3
   Q+1\right)\right.\right.\\
   &\left.\left.+459 (R (6 R+43)-55 Q+6) (R-Q)\right]\right\}\,\\
\end{split}
\end{equation}

\subsection {Surfactant-free (clean) drop with viscosity ratio $\lambda=1$}

\begin{equation}
\begin{split}
F_{20}^{(2)}=&\left[{13660416 (2 R+3)^5
   (4 R+5)}\right]^{-1}125 s^4 \left(1432192 R^6+3288416 R^5-\right.\\
   &\left.8 R^4 (174680 Q+242639)+10 R^3
   (414304 Q-1246403)+\right.\\
  & \left. R^2 \left(-584672 Q^2+16072364 Q-11634545\right)+R \left(-8185408
   Q^2+13647766 Q+113294\right)-\right.\\
   &\left. 9220200 Q^2+3514830 Q+3215505\right)
\end{split}
\end{equation}
\begin{equation}
\begin{split}
F_{40}^{(2)}=&\left[{29597568 (2 R+3)^5 (4 R+5)}\right]^{-1}5 s^4 \left(64383232 R^6+203157536 R^5\right.\\
&\left.+R^4 (251051632-251193984 Q)+R^3
   (61496714-498730704 Q)\right.\\
   &\left. +R^2 \left(230840192 Q^2-242910304 Q-156343191\right)\right.\\
 &  \left. +2 R
   \left(91644014 Q^2+92491121 Q-60629679\right)-15 \left(9376026 Q^2-17391806
   Q+1933249\right)\right)
\end{split}
\end{equation}
\begin{equation}
\begin{split}
F_{60}^{(2)}=&\left[{47811456 (2 R+3)^4 (4 R+5)}\right]^{-1}25 s^4 \left(5710688 R^5+15161920 R^4\right.\\
&\left.+R^3 (17871454-24878432 Q)+R^2
   (9427895-43762456 Q)\right.\\
   &\left. +R \left(24695000 Q^2-24601697 Q+395633\right)+5 \left(3990800
   Q^2+788917 Q-783718\right)\right)
\end{split}
\end{equation}
\begin{equation}
\begin{split}
F_{80}^{(2)}=F_{40}^{(1)} \frac{5 s^2 \left(28772 R^2+31281 R-88825 Q+28772\right)}{339768 (2 R+3)^2}
\end{split}
\end{equation}
In the above expressions $s=2$ corresponds to the linear field $(-x,-y,2z)$

\section*{Acknowledgments}
\label{acknowledgments}
CS gratefully acknowledges support by Comsol Inc. AT has been supported by the Swedish Research Council under Grant 2015-04998 and PV has been supported in part by NSF award  CBET 1704996.\\
We thank Shravan Veerapaneni for the valuable discussions, inputs and suggestions for starting this work and also for sharing the fast code for the rotation of the spherical harmonics expansions (Z. Gimbutas and S. Veerapaneni \cite{gimbutas}) which has been used in the present work.\\
We also would like to thank Herve Nganguia for sharing data from his theory \cite{Nganguia2013} and Shubhadeep Mandal \cite{Mandal_2016} for sharing his code for small deformation theory.

\clearpage
\newpage

\clearpage
\newpage
\section*{Bibliography}
\bibliography{biblio}

\end{document}